\documentclass[a4paper,10pt]{amsart}
\usepackage[margin=1in]{geometry}
\usepackage{enumerate, amsmath, amsfonts, amssymb, amsthm, mathtools, thmtools, wasysym, graphics, graphicx, xcolor, frcursive,xparse,comment,ytableau,stmaryrd,bbm,array,colortbl,tensor,arydshln,leftidx, mathrsfs, caption, subcaption}
\usepackage[all]{xy}
\usepackage[boxed,norelsize]{algorithm2e}
\usepackage{hyperref}
\usepackage{cancel}
\usepackage{calligra}
\usepackage[vcentermath]{youngtab}
\ytableausetup{centertableaux}

\usepackage{bm}
\DeclareMathAlphabet{\mathcalligra}{T1}{calligra}{m}{n}

\makeatletter
\newcommand{\specialcell}[1]{\ifmeasuring@#1\else\omit$\displaystyle#1$\ignorespaces\fi}
\makeatother

\usepackage{url, hypcap}
\hypersetup{colorlinks=true, citecolor=darkblue, linkcolor=darkblue}
\definecolor{join}{RGB}{0,77,178}
\usepackage{tikz}
\usetikzlibrary{calc,through,backgrounds,shapes,matrix,cd}

\definecolor{darkblue}{rgb}{0.0,0,0.7} 
\newcommand{\darkblue}{\color{darkblue}} 
\definecolor{darkred}{rgb}{0.7,0,0} 
 \definecolor{lightgrey}{rgb}{0.7,0.7,0.7}

\definecolor{meet}{RGB}{255,205,111}
\definecolor{join}{RGB}{0,77,178}

\DeclareFontEncoding{OT2}{}{}

\newtheorem{theorem}{Theorem}[section]
\newtheorem{proposition}[theorem]{Proposition}

\newtheorem{lemma}[theorem]{Lemma}

\theoremstyle{definition}
\newtheorem{definition}[theorem]{Definition}
\newtheorem{example}[theorem]{Example}

\newtheorem{conjecture}[theorem]{Conjecture}

\newtheorem{remarkx}[theorem]{Remark}
\newenvironment{remark}
  {\pushQED{\qed}\remarkx}
  {\popQED\endremarkx}

\usepackage[procnames]{listings}
\usepackage[all]{xy}
\usepackage[T1]{fontenc}

\usepackage[cal=boondoxo]{mathalfa}
\usepackage[colorinlistoftodos]{todonotes}
\usepackage{xifthen}
\usepackage{xstring}

\newcommand{\defn}[1]{\emph{\darkblue #1}}

\numberwithin{equation}{subsection}

\usetikzlibrary{math}
\SetKwInput{kwset}{set}
\usetikzlibrary{arrows,backgrounds,calc,trees}
\pgfdeclarelayer{background}
\pgfsetlayers{background,main}
\renewcommand{\mod}{\operatorname{mod}}

\newcommand{\C}{\mathbb{C}}
\newcommand{\Z}{\mathbb{Z}}

\newcommand{\Sn}{\mathfrak{S}}

\newcommand{\A}{\mathcal{A}}

\newcommand{\R}{\mathcal{R}}

\newcommand{\NC}{\mathrm{NC}}

\newcommand{\Cat}{\operatorname{Cat}}
\newcommand{\shuffle}{\operatorname{Sf}}

\renewcommand{\a}{\mathbf{a}}

\newcommand{\Fix}{\mathrm{Fix}}

\newcommand{\CC}{\mathbb{C}}

\newcommand{\GL}{\mathrm{GL}}

\newcommand{\Park}{\operatorname{Park}}
\newcommand{\NN}{\operatorname{NN}}

\DeclareMathAlphabet{\mathcal}{OMS}{cmsy}{m}{n}

\allowdisplaybreaks

\usepackage[nameinlink]{cleveref}
\Crefformat{footnote}{#2\footnotemark[#1]#3}
\Crefname{conjecture}{Conjecture}{Conjectures}
\Crefname{assumption}{Assumption}{Assumptions}
\Crefname{subsection}{Subsection}{Subsections}
\Crefname{remarkx}{Remark}{Remarks}
\title[Parking Spaces for Complex Reflection Groups]{Parking Spaces for Complex Reflection Groups}

\author[Stack]{Jason Stack}
\address[Stack]{University of Texas at Dallas}
\email{jason.stack@utdallas.edu}

\begin{document}
\begin{abstract}
    We answer an open problem of~\cite{armstrong2015parking} and~\cite{rhoades2014parking}, extending their work to irreducible well--generated complex reflection groups $W$. We define a combinatorial $W$-noncrossing parking space and an algebraic $W$-parking space for such $W$, and exhibit a $(W \times C)$-equivariant isomorphism between the two. As a consequence of this isomorphism, we enumerate the $W$-noncrossing parking functions. Finally, we extend our results to the Fuss case. We prove the results for all such complex reflection groups except $G_{34}$, $E_7,$ and $E_8$.
\end{abstract}
\maketitle
\section{Introduction}
\label{sec: intro}
\subsection{Noncrossing Partitions}
\label{subsec:catcomb}
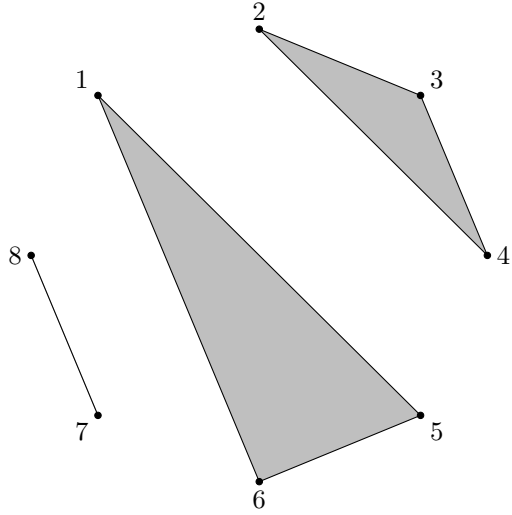
\begin{figure}[h]
    \centering
    \begin{tikzpicture}
%Connect zero block 
\draw[fill=lightgray] (-360/8*1:-3cm)--
(-360/8*5:-3cm) --
(-360/8*6:-3cm)--
(-360/8*1:-3cm);

\draw[fill=lightgray] (-360/8*2:-3cm)--
(-360/8*3:-3cm) --
(-360/8*4:-3cm)--
(-360/8*2:-3cm);

\draw[fill=lightgray] (-360/8*7:-3cm)--
(-360/8*8:-3cm);

\foreach \a/\n in {1/1,2/2,3/3,4/4,5/5,6/6,7/7,8/8}
{
\draw[fill](-360/8*\n:-3cm)circle(1.2pt)
node[anchor=-\n*360/8] (\a) {$\a$};
}

\end{tikzpicture}
    \caption{A noncrossing partition in type $A_{7}$ with 
    $\pi=\{\{1,5,6\},\{2,3,4\},\{7,8\}\}.$}
    \label{fig: Type A noncrossing partition}
\end{figure}

One of the most studied noncrossing objects in Coxeter--Catalan combinatorics are the \defn{noncrossing partitions}. The noncrossing partitions in type $A$ are partitions $\pi=\{B_1,\dots,B_r\}$ of the set $[n]\coloneq\{1,2,\dots,n\}$ such that there are no $a,b,c,d \in [n]$, $a<b<c<d$ where $a,c\in B_i$ and $b,d\in B_j$ for $i\neq j$. As shown in~\Cref{fig: Type A noncrossing partition}, we can visualize this condition by drawing $n$ nodes around a circle. We label the nodes in clockwise order by the elements of $[n]$ with the condition that the intersection of the convex hulls of distinct blocks is empty, hence the name \emph{noncrossing} partitions. 

We can extend this definition to finite Coxeter groups. Let $W$ be a finite Coxeter group with reflection representation $V$ and fix a Coxeter element $c \in W$. The poset of all \defn{$W$-noncrossing partitions}, $\NC(W,c)$, is defined as the interval $[e,c]_T$ in the absolute order (the Cayley graph of $W$ generated by its reflections $T$). The $W$-noncrossing partitions are counted by the \defn{Coxeter-Catalan number}
\[
\Cat(W)\coloneq\prod_{i=1}^n \frac{h+d_i}{d_i},
\]
where $d_1\leq d_2 \leq \dots \leq d_n$ are the degrees of a set of algebraically independent homogeneous polynomials generating the algebra of invariants $\operatorname{Sym}(V^*)^W$ and $h=d_n$ is the Coxeter number of $W$~\cite{reiner1997classical,bessis2003dual}. There is a Fuss generalization of these noncrossing partitions, called \defn{$k$-$W$-noncrossing partitions}, given by multichains in $\NC(W,c)$ of length $k$~\cite{armstrong2009thesis}. Noncrossing partitions and their Fuss generalizations can be further defined for irreducible well--generated complex reflection groups (see~\Cref{def:Well_generated_crg,def:ncpartitions_complex,def:Fuss_ncpartitions_complex}).

As we explain in the next section, the noncrossing partitions are used to define the \defn{$W$-noncrossing parking functions} for real reflection groups $W$~\cite{armstrong2015parking,edelman1980chain}. The object of this paper is to extend these noncrossing parking functions to irreducible well--generated complex reflection groups.

\subsection{Parking Functions}\label{sec:intro_parking}
We first recall classical (nonnesting) parking functions and the more recent noncrossing parking functions for real reflection groups~\cite{edelman1980chain,armstrong2015parking}.
\subsubsection{Nonnesting Parking Functions}
A  \defn{classical parking function} is a sequence of positive integers, $(a_1,\dots,a_n)$, such that their increasing rearrangement, $(b_1,b_2,\dots,b_n)$, satisfies $b_i\leq i$ for all $i=1,\dots,n$. The set of all parking functions, which we denote $\operatorname{Park}_n$, becomes the \defn{classical parking space} when endowed with the action of $\Sn_n$ via $(w\cdot f)(i)=f(w^{-1}(i))$. Note that $\operatorname{Park}_n$ is of size $(n+1)^{n-1}$ and has Catalan many $\Sn_n$-orbits.
\begin{example}
    The table below shows the elements of $\operatorname{Park}_3$ ordered in rows by their $\Sn_3$-orbits:
    \begin{center}
    \begin{tabular}{|c|ccccc|}\hline
        111 &&&&&  \\\hline
         112 &121&211&&& \\\hline
         113 &131&311&&& \\\hline
         122 &212&221&&& \\\hline
         123 &132&312&321&231&213\\\hline
    \end{tabular}.
    \end{center}
    There are $(3+1)^{3-1}=16$ elements and $5=\frac{1}{4}\binom{6}{3}$ orbits.
\end{example}
The action of $\Sn_n$ on $\Park_n$ can be generalized to crystallographic real reflection groups, or \defn{Weyl groups}. Let $W$ be a Weyl group acting irreducibly on $V\cong \mathbb{R}^n$. As $W$ is crystallographic, we may choose the simple roots, positive roots, root system, and root lattice, $\Delta$, $\Phi^+$, $\Phi$, and $Q$, for $W$ with
\[
\Delta \subseteq \Phi^+ \subseteq \Phi \subseteq Q \subseteq V.
\]
The quotient $Q/(h+1)Q$ has cardinality $(h+1)^n$ and carries a $W$-action induced by the action of $W$ on $V$. We call this space the \defn{standard $W$-parking space}~\cite[Section 2.4, Section 7.3]{Haiman1994conjectures}.

There is a more combinatorial definition of the parking space available.  Recall there is a standard partial order on $\Phi^+$ given by $\alpha \leq \beta$ if and only if $\beta - \alpha \in \mathbb{N}\Phi^+$.  Define the set of all \defn{$W$-nonnesting partitions}, $\NN(W)$, to be the set of all antichains in the root poset $(\Phi^+,\leq)$. As with the noncrossing partitions, the $W$-nonnesting partitions are enumerated by $\Cat(W)$~\cite[Theorem 2.4]{athanasiadis1998noncrossnonnest}. Further, the $W$-orbits on $Q/(h+1)Q$ can be parametrized by the $W$-nonnesting partitions (see~\cite[Section 4]{Cellini2002Adnil},~\cite[Section 5]{Sommers2005stable} and~\cite{shi1997signtypes}). Define, for any $\pi\in\NN(W)$, the \defn{parabolic subgroup generated by $\pi$}, $W_\pi\coloneq\langle t\leq \pi: t\in T \rangle$, where $T$ is the set of all reflections in $W$. The space of \defn{$W$-nonnesting parking functions} is the set of cosets
\[
\Park^{\NN} (W)\coloneq\{wW_\pi:w\in W\text{ and } \pi\in \NN(W)\},
\]
under left multiplication by elements of $W$.  By the discussion above, $\Park^{\NN} (W)$ is in $W$-equivariant bijection with $Q/(h+1)Q$~\cite[Lemma 4.1, Theorem 4.2]{athanasiadis2004catw}.

\subsubsection{Noncrossing Parking Functions}
\label{sec:noncrossing_parking_intro}
\begin{figure}[h]
    \centering
\begin{tikzpicture}
%Connect zero block 
\draw[fill=lightgray] (-360/8*1:-3cm)--node[above,xshift=-.25cm,yshift=-.5cm]{\tiny \color{blue} $2,5,6$}
(-360/8*5:-3cm) --
(-360/8*6:-3cm)--
(-360/8*1:-3cm);

\draw[fill=lightgray] (-360/8*2:-3cm)--
(-360/8*3:-3cm) --
(-360/8*4:-3cm)--node[above,rotate=-40]{\tiny \color{blue} $3,4,7$}
(-360/8*2:-3cm);

\draw[fill=lightgray] (-360/8*7:-3cm)--node[above,rotate=-70]{\tiny \color{blue} $1,8$}
(-360/8*8:-3cm);

\foreach \a/\n in {1/1,2/2,3/3,4/4,5/5,6/6,7/7,8/8}
{
\draw[fill](-360/8*\n:-3cm)circle(1.2pt)
node[anchor=-\n*360/8] (\a) {$\a$};
}

\end{tikzpicture}
    \caption{A noncrossing parking function in type $A_{7}$ with noncrossing partition $\pi$ as in~\Cref{fig: Type A noncrossing partition}.}
    \label{fig: Type A noncrossing parking fnct}
\end{figure}
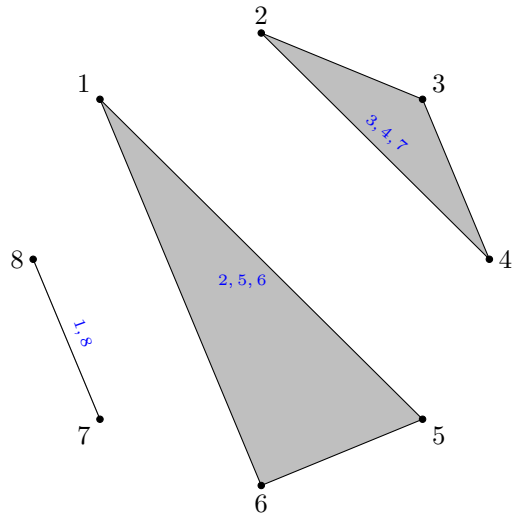

In parallel fashion, for $W$ a Coxeter group with Coxeter element $c$, we may use the $W$-noncrossing partitions from~\Cref{subsec:catcomb} to define the space of \defn{$W$-noncrossing parking functions}~\cite{armstrong2015parking} as the set of cosets
\[
\Park^{\NC} (W)\coloneq\{wW_\pi:w\in W\text{ and } \pi\in \NC(W)\},
\]
again under left multiplication of a coset by elements of $W$. Here, $W_\pi$ is defined analogously for $\pi\in \NC(W)$. The $W$-noncrossing parking functions carry an action of $W\times C$, where $C=\langle c \rangle$. We describe the $W$-noncrossing parking functions in more detail in~\Cref{subsec:ncpf_in_rrg}.
\begin{example}
Introduced by P.~Edelman~\cite{edelman1980chain} as 2-partitions, the set of all $A_{n-1}$-noncrossing parking functions can be viewed as the set of tuples ${(B_1,L_1),(B_2,L_2), \dots,(B_r,L_r)}$ where $B_i,L_i\subset [n]$ such that 
\begin{itemize}
    \item $\pi\coloneq\{B_1,\dots, B_r\}$ \text{ is a \emph{noncrossing} partition of } $[n]$,
    \item $\{L_1,\dots,L_r\}$ \text{ is a set partition of } $[n]$, \text{ and }
    \item $|B_i| = |L_i|$ \text{ for each } $i=1,\dots, r$.
\end{itemize}
We may visualize this by decorating our model in~\Cref{fig: Type A noncrossing partition} with labels $L_i$ associated to each block $B_i\in\pi$, as shown in~\Cref{fig: Type A noncrossing parking fnct}. 
\end{example}

\subsubsection{Fuss Analogues}

There is a Fuss analogue of the classical parking functions dependent on a \defn{Fuss parameter} $k\in \Z_{\geq1}$. A \defn{classical $k$-parking function} is a sequence of positive integers $(a_1,\ldots,a_n)$ whose positive rearrangement $(b_1,\ldots,b_n)$ satisfies $b_i\leq (i-1)k+1$. Note that when $k=1$ we recover $\Park_n$. Denote the set of all $k$-parking functions by $\Park_n^k$. This set again carries an action by $\Sn_n$ via coordinate permutation. 

Generalizing to the crystallographic case, the finite torus $Q/(kh+1)Q$ carries a $W$-action and plays the role of the Fuss parking functions. There is a parametrization of $Q/(kh+1)Q$ by certain geometric multichains in $\NN(W)$ arising from the hyperplane arrangement of $W$~\cite{athanasiadis1999geochains}. This suggests an analogous definition using the noncrossing partitions, which we refer to as the \defn{$k$-$W$-noncrossing parking functions}~\cite{rhoades2014parking}. Similarly to the $k=1$ case, the $k$-$W$-noncrossing parking space carries an action of $W\times \Z_{kh}$. We discuss this space in more detail in~\Cref{sec:Fuss_background}.

\subsection{Algebraic Parking Space}
\label{sec:alg_park_intro}
D.~Armstrong, V.~Reiner, and B.~Rhoades construct a second parking space purely algebraically~\cite{armstrong2015parking,rhoades2014parking}. Let $W$ be a Coxeter group with reflection representation $V$ and fix a Coxeter element $c\in W$. The \defn{algebraic $W$-parking space} (or the \defn{$k$-algebraic $W$-parking space} in the Fuss case) relies on a choice of a degree $p$ \defn{homogeneous system of parameters}, $(\Theta_p)$, and the quotient $S/(\Theta_p - {\bf x})$, where $S\coloneq\operatorname{Sym}{V^*}$. We choose $\Theta_{h+1}$ of degree $h+1$ (or $kh+1$) and have that this quotient carries a $(W\times C)$-action (or $(W\times \Z_{kh})$-action). D.~Armstrong, V.~Reiner, and B.~Rhoades construct isomorphisms between the $W$-noncrossing parking functions and the algebraic $W$-parking space that are invariant under these actions~\cite{armstrong2015parking}. B.~Rhoades does the same in the Fuss case~\cite{rhoades2014parking}.

In this paper, we extend the work of~\cite{armstrong2015parking} and~\cite{rhoades2014parking} to irreducible well--generated complex reflection groups. For these groups, we have well defined notions of Coxeter elements, a Coxeter number, and noncrossing partitions. We are then able to define a purely combinatorial noncrossing parking space in analogy with that defined for finite real reflection groups. We then define a second parking space purely algebraically, extending the construction from the real case~\cite{armstrong2015parking,rhoades2014parking}. We give an example of the existence of specific hsops for the exceptional group $G_{10}$ in~\Cref{ex:G10_hsop}. We exhibit, case--by--case, an isomorphism between the combinatorial and algebraic spaces, and we enumerate the $W$-noncrossing parking functions. We will formally define the algebraic $W$-parking space in~\Cref{sec: defs and thms}.

Let $W$ be an irreducible well--generated complex reflection group not of type $G_{34}$, $E_7$, or $E_8$. Let $h$ be the Coxeter number of $W$ and $n$ its rank. 

\begin{theorem}
For all such $W$, there is a $(W\times C)$-equivariant isomorphism between the $W$-noncrossing parking functions and the algebraic $W$-parking space. 
\end{theorem}

\begin{theorem}
For all such $W$ and $k\geq 1$, there is a $(W\times \Z_{kh})$-equivariant isomorphism between the $k$-$W$-noncrossing parking functions and the $k$-algebraic $W$-parking space. 
\end{theorem}

Combining the results above with computations using~\cite{sagemath,JM15,GH96}, we obtain the following results:

\begin{theorem}
For all such $W$, there are $(h+1)^n$ $W$-noncrossing parking functions.
\end{theorem}

\begin{theorem}
For all such $W$, there are $(kh+1)^n$ $k$-$W$-noncrossing parking functions.
\end{theorem}

\subsection{Outline}
We first introduce the necessary noncrossing objects in real reflection groups and the required background for complex reflection groups in~\Cref{sec: background}. We then present the noncrossing parking functions, both the combinatorial and algebraic definitions, for irreducible well--generated complex reflection groups in~\Cref{sec: defs and thms} along with our main results. We state analogous definitions and results in the Fuss case as well. 

In~\Cref{sec:gd1n} and~\Cref{sec:gddn} we describe an isomorphism between the two parking spaces for the infinite families $G(d,1,n)$ and $G(d,d,n)$, respectively. Similarly, in~\Cref{sec:fuss_gd1n} and~\Cref{sec:fuss_gddn} we describe an isomorphism between the two parking spaces in the Fuss case for the two infinite families $G(d,1,n)$ and $G(d,d,n)$, respectively.

We discuss future projects in~\Cref{sec:future_work}.

\section{Background}
\label{sec: background}
\subsection{Noncrossing Partitions in Real Reflection Groups}
\label{subsec: ncp_in_rrg}

We will begin with a review of noncrossing partitions for real reflection groups. Let $(W,S)$ be a finite Coxeter system, with simple generators $S=\{s_1,\ldots,s_n\}$ and fix a \defn{Coxeter element} $c$, a regular element with order the \defn{Coxeter number} $h$---this definition is discussed further in~\Cref{sec:regular}; for the moment, the reader could imagine that $c$ is just a product of the simple reflections, in some order. Let $T\coloneq\{t_\alpha\}$ denote the set of all reflections belonging to $W$ with corresponding reflecting hyperplanes $H_\alpha$. Let $\operatorname{Cox}(\Phi)\coloneq\{H_\alpha\}_{\alpha \in \Phi}$ be the arrangement of reflecting hyperplanes and $\mathcal{L}$ its lattice of intersecting subspaces, called \defn{flats}, ordered by reverse-inclusion. 
\begin{definition}\label{def:real_absolute_order}
    The \defn{reflection length}, $\ell_T$, of $w\in W$ is the smallest $r$ such that $w=t_1\cdots t_r$ with $t_i\in T$ for all $i$. We define the \defn{absolute order}, $\leq_T$, on $W$ by
    \[
    u \leq_T v \iff \ell_T(v) = \ell_T(u) + \ell_T(u^{-1}v).
    \]
    
\end{definition}
This poset is graded with rank function $\ell_T$. Its unique minimal element is the identity $e\in W$, and among its maximal elements are the Coxeter elements (which form a single conjugacy class if $W$ is crystallographic)~\cite{springer1974regular}. %\nathan{probably should cite Springer ``Regular elements''}

\begin{definition}[\cite{bessis2003dual}\cite{watt2002Kpi}]
\label{def:realnoncrossing-partition}
The poset of \defn{$W$-noncrossing partitions}, $\NC(W,c)$, is the interval $[e,c]_T$ in the absolute order.
\end{definition}
 In~\cite{brady2002partial}, T.~Brady and C.~Watt showed there is an embedding
\[
\NC(W,c)\xhookrightarrow{}\mathcal{L}
\]
given by 
\[
\pi \mapsto \operatorname{Fix}(\pi)\coloneq\{v\in V:\pi v = v\}.
\]  We will drop the dependence on $c$ in $\NC(W,c)$ as there is an isomorphism between $\NC(W,c)$ and $\NC(W,c')$ for any Coxeter elements $c,c'$~\cite[Corollary 1.5]{reiner2017non}. Furthermore, the map $w\mapsto cwc^{-1}$, or the square of the \defn{Kreweras complement} $w \mapsto cw^{-1}$, is a poset isomorphism of $\NC(W,c)$.

\begin{example}
\label{ex:real-nc-partitions}
In the classical (real) types $A$, $B$, and $D$, we can model the noncrossing partitions pictorially~\cite{armstrong2015parking}; these models are illustrated in~\Cref{fig:nc_pics}.  We will occasionally write $\overline{i}\coloneq-i$ for ease of notation.

\begin{itemize}
\item Recall from~\Cref{sec: intro} that the $A_{n-1}$-noncrossing partitions are simply the noncrossing set partitions of $[n]$ .

\item   The Weyl group of type $B_n$, is the \defn{hyperoctahedral group}, $\Sn_n^\pm$, of $n\times n$ signed matrices. We will think of these matrices as permutations of the set $\pm[n]\coloneq\{-n,-n+1,\ldots,-1,1,\ldots,n\}$. The $B_n$-noncrossing partitions are ordinary noncrossing partitions $\pi=\{B_1,\ldots,B_r\}$ of $\pm[n]$ with the added condition that for each $B_i\in \pi$ there exists $j$ such that $B_i = -B_j$. The \defn{zero block} (if it exists) is the unique block, denoted $B_0$, satisfying $B_0=-B_0$. These additional conditions give a central symmetry of the model.
    
\item The Weyl group of type $D_n$ is the subgroup of the hyperoctahedral group consisting of those elements with an even number of $-1$ entries in their matrix representation. The $D_n$-noncrossing partitions are the $B_n$-noncrossing partitions with the additional condition that if there exists a zero block, $B_0 \in \pi$, then $\pm n \in B_0$ and $|B_0|\geq 4$. 
\end{itemize}
We give an example of a noncrossing partition in each type in~\Cref{fig:nc_pics}.
\end{example}
    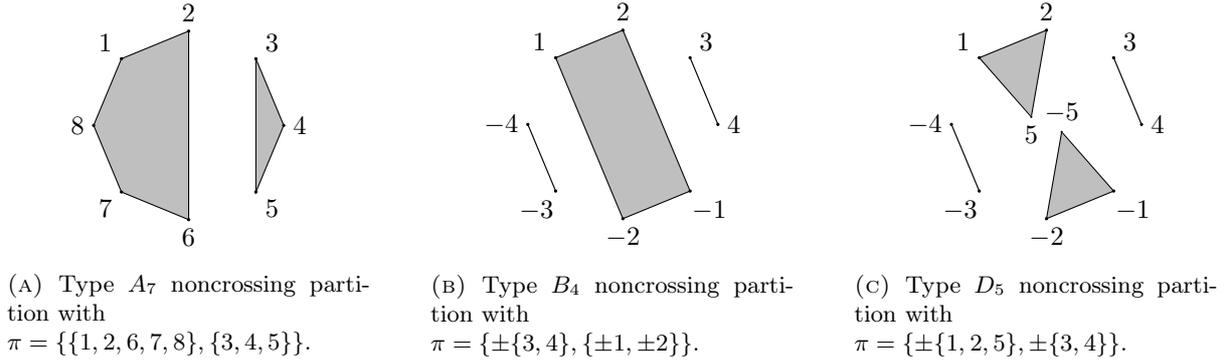
\begin{figure}[h]
        \centering
        \begin{subfigure}[h]{0.3\textwidth}
            \centering
            \begin{tikzpicture}

\draw[fill=lightgray] (-360/8*1:-1.25cm)--
(-360/8*2:-1.25cm) --
(-360/8*6:-1.25cm)--
(-360/8*7:-1.25cm)--
(-360/8*8:-1.25cm)--
(-360/8*1:-1.25cm);

\draw[fill=lightgray] (-360/8*3:-1.25cm)--
(-360/8*4:-1.25cm) --
(-360/8*5:-1.25cm)--
(-360/8*3:-1.25cm);

\foreach \a/\n in {1/1,2/2,3/3,4/4,5/5,6/6,7/7,8/8}
{
\draw[fill](-360/8*\n:-1.25cm)circle(.4pt)
node[anchor=-\n*360/8] (\a) {$\a$};
}

\end{tikzpicture}
            \caption{Type $A_7$ noncrossing partition with\\ $\pi=\{\{1,2,6,7,8\},\{3,4,5\}\}$.}
        \end{subfigure}\hfill
        \begin{subfigure}[h]{0.3\textwidth}
            \centering
            \begin{tikzpicture}
        %Connect zero block 
\draw[fill=lightgray] (-360/8*1:-1.25cm)--
(-360/8*2:-1.25cm) --
(-360/8*5:-1.25cm)--
(-360/8*6:-1.25cm)--
(-360/8*1:-1.25cm);

\draw[fill=lightgray] (-360/8*3:-1.25cm)--
(-360/8*4:-1.25cm);

\draw[fill=lightgray] (-360/8*7:-1.25cm)--
(-360/8*8:-1.25cm);

\foreach \a/\n in {1/1,2/2,3/3,4/4,-1/5,-2/6,-3/7,-4/8}
{
\draw[fill](-360/8*\n:-1.25cm)circle(.4pt)
node[anchor=-\n*360/8] (\a) {$\a$};
}

\end{tikzpicture}
            \caption{Type $B_4$ noncrossing partition with\\ $\pi = \{\pm\{3,4\},\{\pm1,\pm2\}\}$.}
        \end{subfigure}\hfill
        \begin{subfigure}[h]{0.3\textwidth}
            \centering
            \begin{tikzpicture}

\node[circle,draw=black, fill=black, inner sep=0pt,minimum size=.4pt] (5) at (-.2,.1) {};

\node[circle,draw=black, fill=black, inner sep=0pt,minimum size=.4pt] (-5) at (.2,-.1) {};

\draw[fill=lightgray] (-360/8*1:-1.25cm)--
(-360/8*2:-1.25cm) --
(5)node[below]{$5$}--
(-360/8*1:-1.25cm);

\draw[fill=lightgray] (-360/8*5:-1.25cm)--
(-360/8*6:-1.25cm) --
(-5)node[above]{$-5$}--
(-360/8*5:-1.25cm);

\draw[fill=lightgray] (-360/8*3:-1.25cm)--
(-360/8*4:-1.25cm);

\draw[fill=lightgray] (-360/8*7:-1.25cm)--
(-360/8*8:-1.25cm);

\foreach \a/\n in {1/1,2/2,3/3,4/4,-1/5,-2/6,-3/7,-4/8}
{
\draw[fill](-360/8*\n:-1.25cm)circle(.4pt)
node[anchor=-\n*360/8] (\a) {$\a$};
}

\end{tikzpicture}
            \caption{Type $D_5$ noncrossing partition with \\$\pi=\{\pm\{1,2,5\},\pm\{3,4\}\}$.}
        \end{subfigure}
        \caption{Noncrossing partitions in real reflection groups for types $A$, $B$, and $D$.}
        \label{fig:nc_pics}
    \end{figure}

\subsection{Noncrossing Parking Functions in Real Reflection Groups}\label{subsec:ncpf_in_rrg}
We follow the construction of the noncrossing parking functions from~\cite{armstrong2015parking}. As in type $A$ from~\Cref{sec:noncrossing_parking_intro}, we decorate the $W$-noncrossing partitions to produce the $W$-noncrossing parking functions.

\begin{definition}
    The set of \defn{$W$-noncrossing parking functions} is defined as
    \[
    \Park^\NC(W)\coloneq\{w W_\pi : w \in W \text{ and } \pi \in \NC(W)\}.
    \]
\end{definition}
We will use the notation $[w,\pi]$ to denote the element $wW_\pi \in \Park^\NC(W)$. This set carries a $(W\times C)$-action given by
\[
(v,c^p)[w,\pi] = [vwc^{-p},c^p\pi c^{-p}],
\]
for any $p\in \mathbb{N}$ and any $v \in W$.

Equivalently, we may define the $W$-noncrossing parking functions in the language of elements in $W\times \NC(W)$ or $W\times \mathcal{L}$ via the embedding described in~\Cref{def:realnoncrossing-partition} which agrees with D.~Armstrong, V.~Reiner, and B.~Rhoades' definition~\cite[Definition 2.6]{armstrong2015parking}. 
\begin{example}
    In the classical real types, we can pictorially represent the noncrossing parking functions using the corresponding visualization of noncrossing partitions from~\Cref{ex:real-nc-partitions}~\cite{armstrong2015parking}; examples are illustrated in~\Cref{fig:parkingpics}. Consider $[w,\pi]\in W$ for $W$ of type $A$, $B$, or $D$. We have a noncrossing partition $\pi$ and an element $w\in W$ only being considered up to its coset $wW_\pi$.  Thus, we can picture $w$ as a function on the blocks of $\pi$.  Taking the noncrossing models from~\Cref{fig:nc_pics}, we may represent a $W$-noncrossing parking function $[w,\pi]\in\Park^\NC (W)$ by decorating each block with its image under $w$.
    
    \end{example}
    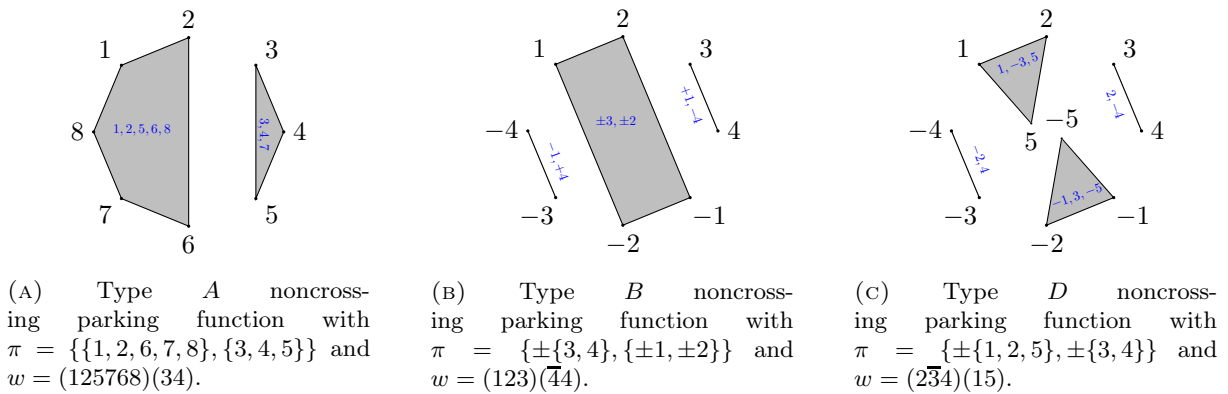
\begin{figure}[h] \label{fig:real_parking_function_examples}
        \centering
        \begin{subfigure}[h]{0.3\textwidth}
            \centering
            \begin{tikzpicture}
\draw[fill=lightgray] (-360/8*1:-1.25cm)--node[below,xshift=-.175cm,yshift=-.85cm]{\scalebox{.5}{\color{blue} $1,2,5,6,8$}}
(-360/8*2:-1.25cm) --
(-360/8*6:-1.25cm)--
(-360/8*7:-1.25cm)--
(-360/8*8:-1.25cm)--
(-360/8*1:-1.25cm);

\draw[fill=lightgray] (-360/8*3:-1.25cm)--
(-360/8*4:-1.25cm) --
(-360/8*5:-1.25cm)--node[above,rotate=-90,yshift=-.1cm]{\scalebox{.5}{ \color{blue} $3,4,7$}}
(-360/8*3:-1.25cm);

\foreach \a/\n in {1/1,2/2,3/3,4/4,5/5,6/6,7/7,8/8}
{
\draw[fill](-360/8*\n:-1.25cm)circle(.4pt)
node[anchor=-\n*360/8] (\a) {$\a$};
}

\end{tikzpicture}
            \caption{Type $A$ noncrossing parking function with $\pi = \{\{1,2,6,7,8\},\{3,4,5\}\}$ and $w = (125768)(34)$.}
        \end{subfigure}\hfill
        \begin{subfigure}[h]{0.3\textwidth}
            \centering
            \begin{tikzpicture}
        %Connect zero block 
\draw[fill=lightgray] (-360/8*1:-1.25cm)--node[below,xshift=.35cm,yshift=-.75cm]{\scalebox{.5}{\color{blue} $\pm 3, \pm 2$}}
(-360/8*2:-1.25cm) --
(-360/8*5:-1.25cm)--
(-360/8*6:-1.25cm)--
(-360/8*1:-1.25cm);

\draw[fill=lightgray] (-360/8*3:-1.25cm)--node[below,rotate=-70]{\scalebox{.5}{ \color{blue} $+1,-4$}}
(-360/8*4:-1.25cm);

\draw[fill=lightgray] (-360/8*7:-1.25cm)--node[above,rotate=-70]{\scalebox{.5}{ \color{blue} $-1,+4$}}
(-360/8*8:-1.25cm);

\foreach \a/\n in {1/1,2/2,3/3,4/4,-1/5,-2/6,-3/7,-4/8}
{
\draw[fill](-360/8*\n:-1.25cm)circle(.4pt)
node[anchor=-\n*360/8] (\a) {$\a$};
}

\end{tikzpicture}
            \caption{Type $B$ noncrossing parking function with $\pi = \{\pm\{3,4\},\{\pm1,\pm2\}\}$ and $w=(123)(\overline{4}4)$.}
        \end{subfigure}\hfill
        \begin{subfigure}[h]{0.3\textwidth}
            \centering
            \begin{tikzpicture}

\node[circle,draw=black, fill=black, inner sep=0pt,minimum size=.4pt] (5) at (-.2,.1) {};

\node[circle,draw=black, fill=black, inner sep=0pt,minimum size=.4pt] (-5) at (.2,-.1) {};

\draw[fill=lightgray] (-360/8*1:-1.25cm)--node[below,rotate=20]{\scalebox{.5}{\color{blue} $1,-3,5$}}
(-360/8*2:-1.25cm) --
(5)node[below]{$5$}--
(-360/8*1:-1.25cm);

\draw[fill=lightgray] (-360/8*5:-1.25cm)--node[above,rotate=20]{\scalebox{.5}{ \color{blue} $-1,3,-5$}}
(-360/8*6:-1.25cm) --
(-5)node[above]{$-5$}--
(-360/8*5:-1.25cm);

\draw[fill=lightgray] (-360/8*3:-1.25cm)--node[below,rotate=-70]{\scalebox{.5}{ \color{blue} $2,-4$}}
(-360/8*4:-1.25cm);

\draw[fill=lightgray] (-360/8*7:-1.25cm)--node[above,rotate=-70]{\scalebox{.5}{\color{blue} $-2,4$}}
(-360/8*8:-1.25cm);

\foreach \a/\n in {1/1,2/2,3/3,4/4,-1/5,-2/6,-3/7,-4/8}
{
\draw[fill](-360/8*\n:-1.25cm)circle(.4pt)
node[anchor=-\n*360/8] (\a) {$\a$};
}

\end{tikzpicture}
            \caption{Type $D$ noncrossing parking function with $\pi=\{\pm\{1,2,5\},\pm\{3,4\}\}$ and $w=(2\overline{3}4)(15)$.}
        \end{subfigure}
        \caption{Parking functions in real reflection groups for types $A$, $B$, and $D$.}
        \label{fig:parkingpics}
    \end{figure}

\subsection{Fuss Analogues}\label{sec:Fuss_background} 
    In~\cite{rhoades2014parking}, B.~Rhoades generalized the noncrossing parking functions for finite real reflection groups to accommodate an extra (integral) \defn{Fuss parameter} $k \geq 1$.  Recall that a \defn{$k$-multichain} in a poset $(P,\leq)$ is a sequence of elements, $(x_1 \leq \cdots \leq x_k)$, in $P$. 
    \begin{definition}[\cite{armstrong2009thesis}]
        Let $W$ be a finite real reflection group. The set $\NC^k(W)$ of \defn{$k$-$W$-noncrossing partitions} consists of all $k$-multichains $(\pi_1\leq_T\cdots\leq_T \pi_k)$ in the poset $\NC(W)$.
    \end{definition}
    When referring to elements of $\NC^k(W)$, the subscript $T$ will be dropped for notational convenience. There is an action of $\Z_{kh}$ in this space which relies on the set of $\ell_T$-additive factorizations of $c$. We say $(\pi_0,\ldots,\pi_k)\in W^{k+1}$ is an \defn{$\ell_T$-additive factorization of $c$} if
    \[
    \pi_0\pi_1\cdots \pi_k=c \text{\hspace{1cm} and \hspace{1cm}} \sum_{i=0}^k\ell_T(\pi_i)=\ell_T(c).
    \]
    
    Let $\NC_k(W)$ be the set of all $\ell_T$-additive factorizations of $c$ of length $k+1$. The following maps due to~\cite{armstrong2009thesis} are mutually inverse bijections between $\NC^k(W)$ and $\NC_k(W)$:
  \begin{align*}
&\partial: &NC^k(W) &\rightarrow NC_k(W) \\
& &\partial(\pi_1 \leq \dots \leq \pi_k) &= 
(\pi_1, \pi_1^{-1} \pi_2, \dots, \pi_{k-1}^{-1} \pi_k, \pi_k^{-1} c)\\
&\int: &NC_k(W) &\rightarrow NC^k(W) \\
& &\int(\pi_0, \pi_1, \ldots, \pi_k) &= (\pi_0 \leq \pi_0 \pi_1 \leq \cdots \leq \pi_0 \pi_1 \cdots \pi_{k-1}).
\end{align*}

Analogously to the square of the Kreweras complement on $\NC(W)$ given by conjugation by $c$, there is an action of $\Z_{kh}$ on $\NC_k(W)$ given by
\[
g\cdot(\pi_0,\ldots,\pi_k)=((c\pi_kc^{-1})\pi_0(c\pi_kc^{-1})^{-1},c\pi_kc^{-1},\pi_1,\pi_2,\ldots,\pi_{k-1}),
\]
with $g$ a distinguished generator of $Z_{kh}$~\cite{armstrong2009thesis}. The maps $\partial$ and $\int$ translate this action into an action on $\NC^k(W)$.

\begin{definition}
    The set of \defn{$k$-$W$-noncrossing parking functions} is defined as 
    \[
    \Park_\NC^k(W)\coloneq\{wW_{\pi_1}\subseteq wW_{\pi_2}\subseteq\cdots\subseteq wW_{\pi_k}:w\in W \text{ and }\pi_i\in\NC(W)\}.
    \]
\end{definition}

Similarly to the $k=1$ case, we will use $[w,(\pi_1\leq\cdots\leq\pi_k)]$ to denote an element of $\Park^k_\NC (W)$. The set $\Park_\NC^k(W)$ carries an action by $W\times \Z_{kh}$, as follows. The action of $W$ is given by
    \[
    v\cdot[w,(\pi_1\leq\cdots\leq\pi_k)] = [vw,(\pi_1\leq\cdots\leq\pi_k)]
    \]
    
   and the $\mathbb{Z}_{kh}$-action on $\Park_\NC^k(W)$ is
    \[
    g\cdot[w,(\pi_1\leq\cdots\leq\pi_k)] = [w\pi_kc^{-1},g\cdot(\pi_1\leq\cdots\leq\pi_k)],
    \]
    where $g\cdot(\pi_1\leq\cdots\leq\pi_k)$ is the action of $g$ on $\NC^k(W)$~\cite{rhoades2014parking}. 
        
In~\cite{rhoades2014parking}, B.~Rhoades describes a unified combinatorial model in types $B$ and $D$, which we extend to the two infinite families of irreducible well--generated complex reflection groups in \Cref{sec:fuss_gd1n} and \Cref{sec:fuss_gddn}, respectively.

\subsection{Complex Reflection Groups}
In this section we will recall the necessary background concerning complex reflection groups, much of which is found in \cite{lehrer2009unitary}.

\subsubsection{Notation} Let $V$ be a finite dimensional complex vector space. A linear transformation $g\in \GL(V)$ is called a \defn{reflection} if the order of $g$ is finite and the subspace $\Fix(g)$ is of codimension 1. When $g$ is a reflection, we call $\Fix(g)$ the \defn{reflecting hyperplane} of $g$.
\begin{definition}
    A \defn{(finite) complex reflection group} is a finite subgroup of $\GL(V)$ that is generated by reflections.
\end{definition}
   For a complex reflection group $W$, we will denote the reflections of $W$ by $\R$ and their associated arrangement of hyperplanes by $\A$. 
\subsubsection{Classification} 
A complex reflection group $W$ is \defn{irreducible} if $V$ is an irreducible $W$-module. The finite irreducible complex reflection groups are classified---there is one infinite family, denoted $G(d,e,n)$ for $e | d$, and 34 exceptional groups, numbered $G_4,\ldots,G_{37}$~\cite{shephard1954finite}.

Since we will be checking our constructions case-by-case, we recall the definition of $G(d,e,n)$ here. 
Let $\mu_d$ be the cyclic group of $\C^\times$ containing $d^{th}$ roots of unity. Consider $G$ a group acting on the set $[n]$ and $V$ a complex inner product space of dimension $n$ having orthonormal basis $e_1,\ldots,e_n$. Then the \defn{standard monomial representation} of $\mu_d\wr G$ is 
\[
(h,g)e_i = h_{g(i)}e_{g(i)}.
\]
With $e$ a divisor of $d$, define
\[
A(d,e,n)\coloneq\{(h_1,\ldots,h_n)\in\mu^n_d:(h_1\cdots h_n)^{d/e}=1\}.
\]
Then, $G(d,e,n)\coloneq A(d,e,n)\rtimes \Sn_n$. Note that $G(d,e,n)$ is a subgroup of $\mu_d\wr\Sn_n$ and as such it may be represented as a group of linear transformations in the standard monomial representation.

\subsubsection{Invariant Theory} 
Let $W$ be an irreducible complex reflection group with reflection representation $V$ of dimension $n$. Let $S\coloneq\operatorname{Sym}(V^*)\cong \C[x_1,\ldots,x_n]$ be the symmetric algebra of $V^*$. Then, there is an action of $W$ on $S$ via $(w\cdot f)(v)=f(w^{-1}v)$ with $f\in S$, $w\in W$, and $v\in V$. The \defn{algebra of invariants}, $S^W$, is generated by algebraically independent homogeneous polynomials~\cite{chevalley1955invariants}. 
\begin{definition}
        For $W$ an irreducible complex reflection group, the \defn{degrees}, $d_1\leq\cdots\leq d_n$, of $W$ are the degrees of the homogeneous polynomials generating $S^W$.
\end{definition}

Let $S^W_+$ be the ideal of $S$ generated by the positive degree elements of $S^W$. Then, the action of $W$ on the \defn{coinvariant algebra}, $S/S^W_+$, is the regular representation. So for any irreducible representation $M$ of $W$, with $M$ having dimension $r$, $S/S_+^W$ contains exactly $r$ copies of $M$.

\begin{definition}
    For $M$ as above, the \defn{exponents} of $M$ are defined to be the degrees, $e_1(M)\leq\cdots\leq e_r(M)$, of the homogeneous components of $S/S_+^W$ containing a copy of $M$. The degrees of $W$ satisfy $d_i=e_i(V)+1$ for $i=1,\ldots,n$. The \defn{codegrees}, $d_1^*\leq\cdots\leq d_n^*$, of $W$ are defined by $d_i^*\coloneq e_{n-i+1}(V^*)-1$.
\end{definition}

\subsubsection{Well--Generated Complex Reflection Groups and Regular Elements}
\label{sec:regular}
To define the noncrossing partitions for complex reflection groups, we restrict to the irreducible \defn{well--generated} complex reflection groups, for which we have the notion of a Coxeter element.

\begin{definition}
\label{def:Well_generated_crg}
    Let $W$ be an irreducible complex reflection group with reflection representation $V$ of dimension $n$. Then, $W$ is \defn{well--generated} if $W$ is generated by $n$ reflections. Equivalently, the degrees and codegrees satisfy $d_i+d_i^*=d_n$. For an irreducible well--generated complex reflection group, we define the \defn{Coxeter number} $h\coloneq d_n$, the highest degree. 
\end{definition}

    The irreducible well--generated complex reflection groups are the two infinite subfamilies---types $G(d,1,n)$ and $G(d,d,n)$---as well as the exceptional groups
    \[
G_i \text{ where } i\in\{4,5,6,8,9,10,14,16,17,18,20,21,23,\dots,30,32,\dots,37\}.
\]

For irreducible well--generated complex reflection groups (where, in contrast to real reflection groups, we have no simple reflections), some care must be taken to define Coxeter elements~\cite{springer1974regular}.  
    \begin{definition}  A vector $v\in V$ is \defn{regular} if it is not contained in any reflection hyperplane $H\in \A$.
       For $\zeta$ a root of unity, an element $c\in W$ is \defn{$\zeta$-regular} if it has a regular eigenvector with eigenvalue $\zeta$. In this case, we say that the order $r$ of $\zeta$ is a \defn{regular number} for $W$.
    \end{definition}

\begin{theorem}[\cite{lehrer1999reflection,lehrer2003invariant}]
    A positive integer $r$ is a regular number for $W$ if and only if $r$ divides the same number of degrees and codegrees. 
\end{theorem}

For $W$ an irreducible well--generated complex reflection group, the Coxeter number $h$ is a regular number for $W$---in this case, the degrees are distinct, $d_i^*=d_{n-i+1}$, and $h=d_n$ is the largest degree, so $h$ only divides $d_n$ and $d_1^*$. It follows that there exists a $\zeta$-regular element for every $h^{\text{th}}$ root of unity $\zeta$. 

\begin{definition}\label{def:coxeterelement}
  The  \defn{Coxeter elements} of $W$ are the $\zeta$-regular elements of $W$, for $\zeta$ an $h^{th}$ root of unity.
\end{definition}
For Weyl groups,~\Cref{def:coxeterelement} is slightly more general than the usual definition of Coxeter element as a product of the simple reflections in any order---for example, any long cycle in the symmetric group $\mathfrak{S}_n$ is a regular element for an $n^{th}$ root of unity---but this just extends the usual definition by conjugation.  For real reflection groups,~\Cref{def:coxeterelement} is even more general, as we can have non-conjugate Coxeter elements~\cite{reiner2017non}.

\section{Results}
\label{sec: defs and thms}

We now define the noncrossing partitions and noncrossing parking functions, both combinatorial and algebraic, for irreducible well--generated complex reflection groups. We then present our main results.
\subsection{Definitions of Noncrossing Objects}

Let $W$ be an irreducible well--generated complex reflection group and let $\R$ be the set of all reflections in $W$. For any $w \in W$, the \defn{reflection length}, $\ell_\R$, of $w$ is the minimal number $r$ such that $w = t_1\cdots t_r$ where $t_i \in \R$. The absolute order, noncrossing partitions, and noncrossing parking functions are defined as they were for real reflection groups in~\Cref{subsec: ncp_in_rrg} and~\Cref{subsec:ncpf_in_rrg}.

\begin{definition}
    For $W$ as above, we define the \defn{absolute order} on $W$ by 
    \[
    u\leq_\R v \iff \ell_\R(v) = \ell_\R(u)+\ell_\R(u^{-1}v).
    \]
\end{definition}
For notational convenience, we drop the $\R$ in $\leq_\R$. Now fix a Coxeter element $c\in W$, as in~\Cref{sec:regular}. 
\begin{definition}\label{def:ncpartitions_complex}
    For $W$ an irreducible well--generated complex reflection group, define the poset $\NC(W,c)$ of \defn{$W$-noncrossing partitions} as the interval $[e,c]_\R$ in the absolute order. 
\end{definition}

Let $\operatorname{Aut}_\R(W)$ be the group of reflection automorphisms of $W$. V.~Reiner, V.~Ripoll, and C.~Stump noted that the action of $\operatorname{Aut}_\R(W)$ preserves the set of Coxeter elements and is transitive on it~\cite[Proposition 1.4]{reiner2017non}. Due to this fact, for any irreducible well--generated complex reflection group $W$, and any two Coxeter elements $c,c'\in W$, there exists a reflection automorphism $\phi$ mapping $c$ to $c'$. As such we have the following result:
\begin{proposition}[\cite{reiner2017non}, Corollary 1.5]
    Let $W$ be an irreducible well--generated complex reflection group, and let $c$ and $c'$ be Coxeter elements. Then $\NC(W,c)$ is isomorphic to $\NC(W,c')$.
\end{proposition}
Consequently, we will drop the $c$ in the notation of $\NC(W,c)$. D.~Bessis proves that the noncrossing partition posets are always lattices~\cite[Lemma 8.6]{bessis2015K}. 

D.~Bessis and R.~Corran give models for the noncrossing partitions of types $G(d,1,n)$ and $G(d,d,n)$, similar to the real types in~\Cref{fig:nc_pics}.   We postpone the details of these models (and their corresponding parking functions) to \Cref{sec: d1n model} and \Cref{sec:gddn_model}.

Having defined the noncrossing partitions in this setting, we may follow the same procedure as in the case of real reflection groups to define the noncrossing parking functions. For $W$ an irreducible well--generated complex reflection group and $\pi\in\NC(W)$, let $W_\pi\coloneq\left<t\leq \pi : t\in \R\right>$ be the \defn{parabolic subgroup of $W$ generated by $\pi$}.

\begin{definition}
\label{def:W_noncrossing_parkingfunction}
    For $W$ an irreducible well--generated complex reflection group, the set of \defn{$W$-noncrossing parking functions} is defined as
    \[
    \Park^\NC (W)\coloneq\{w W_\pi : w \in W \text{ and } \pi \in \NC(W)\}.
    \]
\end{definition}

Similarly to the real case, we use $[w,\pi]$ to denote an element of $\Park^\NC (W)$. There is a $(W\times C)$-action on $\Park^\NC (W)$ given by
\[
(v,c^p)\cdot[w,\pi] = [vwc^{-p},c^p\pi c^{-p}],
\]
for any $p\in\mathbb{N}$ and $v\in W$. By definition of $\Park^\NC (W)$, we have an evident $W$-equivariant bijection
\[
\Park^\NC (W) \cong \bigoplus_{\pi\in\NC(W)}W/W_\pi.
\]

We similarly define the Fuss noncrossing parking space built out of multichains in $\NC(W)$. Let $k\geq 1$.

\begin{definition}\label{def:Fuss_ncpartitions_complex}
    For $W$ an irreducible well--generated complex reflection group, define the set of \defn{$k$-$W$-noncrossing partitions} by
    \[
    \NC^k(W)\coloneq\{(\pi_1\leq\cdots\leq\pi_k):\pi_i\in\NC(W)\}.
    \]
\end{definition}
\begin{definition}
\label{def:Fuss_noncrossing_parkingfunction}
    For $W$ an irreducible well--generated complex reflection group, the set of \defn{$k$-$W$-noncrossing parking functions} is defined as
    \[
    \Park^k_\NC(W)\coloneq\{w W_{\pi_1}\subseteq w W_{\pi_2}\subseteq \cdots \subseteq w W_{\pi_k} : w \in W \text{ and } \pi_i \in \NC(W)\}.
    \]
\end{definition}
We will use $[w,(\pi_1\leq\cdots\leq\pi_k)]$ to denote an element of $\Park^k_\NC (W)$. Similarly to the real case, $\Park_\NC^k(W)$ carries an action by $W\times \Z_{kh}$. The action of $W$ is straightforward: for $v\in W$, we have
\[
v\cdot[w,(\pi_1\leq\cdots\leq\pi_k)]=[vw,(\pi_1\leq\cdots\leq\pi_k)].
\]

To describe the action of $\Z_{kh}$, we first introduce the following set:
\begin{definition}[\cite{KM2013}]
   The \defn{$k$-divisible $W$-noncrossing partitions} are the set
\[
\NC_k(W,c)\coloneq\{(w_0;w_1,\ldots,w_k):w_0w_1\cdots w_k = c \text{ and } \ell_\R(w_0)+\ell_\R(w_1)+\cdots+\ell_\R(w_k)=\ell_\R(c)\}.
\]
A partial order on $\NC_k(W,c)$ is given by
\[
(w_0;w_1,\ldots,w_k)\leq (u_0;u_1,\ldots,u_k) \iff w_i\leq_\R u_i \text{ for all } 1\leq i \leq k. 
\]
\end{definition}

Again, by~\cite{reiner2017non}, we drop the dependence on $c$. When $k=1$ we recover $\NC(W)$. Just as in the real case, the maps $\delta$ and $\int$ give a bijection between $\NC_k(W)$ and $\NC^k(W)$. In~\cite{KM2013}, C.~Krattenhaler and T.~W.~M\"uller give an action of order $kh$ on $\NC_k(W)$ generated by the map 
\[
(w_0,\ldots,w_k)\mapsto ((cw_kc^{-1})w_0(cw_kc^{-1})^{-1},cw_kc^{-1},w_1,w_2,\ldots,w_{k-1}).
\]
Again, $\delta$ and $\int$ again translate this into an action on $\NC^k(W)$. Using this action, we define the action of $\Z_{kh}$ on $\Park^k_\NC (W)$ by
\[
g\cdot[w,(\pi_1\leq\cdots\leq\pi_k)] = [w\pi_kc^{-1},g\cdot(\pi_1\leq\cdots\leq\pi_k)],
\]
where $g\cdot(\pi_1\leq\cdots\leq\pi_k)$ is the $\Z_{kh}$ action on $\NC^k(W)$.

We have thus defined the combinatorial parking spaces, and now turn to the algebraic parking space.

\subsection{Algebraic Parking Space}
\label{sec:alg_park}
The construction of the algebraic parking space in real and complex type relies on the choice of a \defn{homogeneous system of parameters (hsop)}. Let $W$ be a complex reflection group with reflection representation $V\cong \C^n$ and let $S\coloneq\operatorname{Sym}(V^*)$ the algebra of polynomial functions on $V$.
\begin{definition}[\cite{bessis2011cyclic}, Definition 4.1]
    Let $U$ be an $n$-dimensional complex representation of $W$. We say that a collection $\Theta_p=(\theta_1,\ldots,\theta_n)$, with each $\theta_i\in S$ all homogeneous of degree $p$, form a \defn{homogeneous system of parameters (hsop) carrying $U(-p)$} if
    \begin{itemize}
        \item $(\theta_1,\ldots,\theta_n)$ are algebraically independent and the quotient $S/(\Theta_p)$ is finite dimensional over $\C$, and
        \item the span $\C\theta_1+\cdots+\C\theta_n$ is a $W$-stable subspace of the $p^{th}$ homogeneous component of $S$ and carries a $W$-representation equivalent to $U$.
    \end{itemize}
    \end{definition}
\subsubsection{Algebraic Parking Spaces for Real Reflection Groups}
Let $W$ be a real reflection group with Coxeter number $h$ and reflection representation $V$ (considered over $\C$). Let $W$ act on the polynomial algebra $S$ and let $x_1,\ldots,x_n$ be a $\C$-basis for $V^*$.
In~\cite[Section 12]{armstrong2015parking}, D.~Armstrong, V.~Reiner, and B.~Rhoades note a proof due to P.~Etingof, that as a consequence of the representation theory of \defn{rational Cherednik algebras} (see~\cite{etingof2003cherednik}\cite{etingof2012supports}\cite{etingof2010notes}\cite{gordon2003quotient}), there exists a degree $h+1$ hsop $\left(\Theta_{h+1}\right) = (\theta_1,\ldots,\theta_n)$ carrying $V^*(-h-1)$ such that the map $V^*\to S$ given by $x_i\to\theta_i$ is $W$-equivariant.  

\begin{definition}[\cite{armstrong2015parking}, Definition 2.10]
    Let $W$ be an irreducible real reflection group with $\left(\Theta_{h+1}\right) = (\theta_1,\ldots,\theta_n)$ and ${\bf x} = (x_1,\ldots,x_n)$ as above. Then, we define the \defn{algebraic $W$-parking space} by
    \[
    \Park^{alg}(W,\Theta_{h+1})\coloneq S/(\Theta_{h+1}-{\bf x}).
    \]
\end{definition}

This quotient ring has a structure of a $(W\times C)$-representation as $(\Theta_{h+1}-\textbf{x})$ is stable under the action of $W$ by linear substitutions and of $C$ by scalar substitutions $c^p(x_i) = \zeta^{-p}(x_i)$ with $\zeta = e^\frac{2\pi\sqrt{-1}}{h}$. 
\begin{example}
In types $B$ and $D$, the construction of such an hsop is straightforward. In each type, $\Theta_{h+1}=(x_1^{h+1},\ldots,x_n^{h+1})$ gives an hsop carrying $V^*(-h-1)$~\cite{armstrong2015parking}. We will see a similar hsop construction for $G(d,1,n)$ and $G(d,d,n)$.

\end{example}
We may drop the use of $\Theta_{h+1}$ in the notation for the algebraic parking space---writing just $\Park^{alg}(W)$---due to the following proposition:
\begin{proposition}[\cite{armstrong2015parking}, Proposition 2.11]\label{prop:alg_park_iso_btw_hsops}
    For every irreducible real reflection group $W$, and for any choice of hsop $\Theta_{h+1}$, one has an isomorphism of $W\times C$ representations
    \[
    \Park^{alg}(W,\Theta_{h+1})\coloneq S/(\Theta_{h+1}-{\bf x}) \cong S/(\Theta_{h+1}).
    \]
\end{proposition} 

B.~Rhoades extends this construction to the Fuss case. Taking $W$ as before, let $(\Theta_{kh+1})=(\theta_1,\ldots,\theta_n)$ be an hsop carrying $V^*(-kh-1)$ so that $x_i\mapsto \theta_i$ is $W$-equivariant~\cite{rhoades2014parking}.
\begin{definition}[\cite{rhoades2014parking}, Definition 3.4]
    Let $W$ be an irreducible real reflection group with hsop $(\Theta_{kh+1}) = (\theta_1,\ldots,\theta_n)$ as above. We define the \defn{$k$-algebraic $W$-parking space} by
    \[
    \Park_{alg}^k(W,\Theta_{kh+1})\coloneq S/(\Theta_{kh+1}-{\bf x}).
    \]
\end{definition}
Similarly to the $k=1$ case, $(\Theta_{kh+1}-{\bf x})$ is stable under the action of $W \times \Z_{kh}$ by linear substitution and scalar substitution. Using the following proposition, we may drop the dependence on $\Theta_{kh+1}$, as with the $k=1$ case.
\begin{proposition}[\cite{rhoades2014parking}, Proposition 3.5]\label{prop:k_alg_park_iso_btw_hsops}
    For every irreducible real reflection group $W$, and for any choice of hsop $(\Theta_{kh+1})$, one has an isomorphism of $(W\times \Z_{kh})$-representations
    \[
    \Park^{alg}(W,\Theta_{kh+1})\coloneq S/(\Theta_{kh+1}-{\bf x})\cong S/(\Theta_{kh+1}).
    \]
\end{proposition}

\subsubsection{Algebraic Parking Space for Complex Reflection Groups}
We now extend the algebraic parking space to irreducible well--generated complex reflection groups. Let $W$ be an irreducible well--generated complex reflection group with Coxeter number $h$ and reflection representation $V$. Let $S\coloneq \operatorname{Sym}(V^*)$ be the algebra of polynomial functions on $V$. 

We first have a result regarding the existence of hsops for well--generated complex reflection groups. The proof follows from the work of I.~Gordan and S.~Griffeth (see~\cite{gordon2012catalan}, Theorem 2.11 and Subsection 2.18) combined with the freeness theorem for finite Hecke algebras~\cite{etingof2017BMRconj}\footnote{The author would like to thank Stephen Griffeth for pointing out this result and Vic Reiner for suggesting to reach out to Stephen.}. Let $(V^*)^{\sigma_p}$ is the Galois twist of $V^*$ corresponding to $p$---that is, an automorphism of $\C$ mapping $e^\frac{2\pi \sqrt{-1}}{h}$ to $e^\frac{2\pi p\sqrt{-1}}{h}$.
\begin{theorem}[\cite{gordon2012catalan},\cite{etingof2017BMRconj}]\label{thm:hsop_complex}
    If $\gcd(p,h)=1$, then $S$ contains a degree $p$ hsop $\Theta_p$ carrying $(V^*)^{\sigma_p}$. 
\end{theorem}

\begin{remark}
\label{rmk:hsop_inf_fams}
D.~Bessis and V.~Reiner show for the infinite family $G(d,e,n)$ with Coxeter number defined to be the largest degree $d_n=\max{\{(n-1)d,\frac{de}{n}\}}$, that $\Theta_p=(x_1^p,\ldots,x_n^p)$ gives an hsop carrying $V^*(-p)$ if $p\equiv 1 \mod d$~\cite{bessis2011cyclic}. 
\end{remark}

\begin{example}\label{ex:G10_hsop}
    We give an extended example constructing hsops in the exceptional rank 2 complex reflection group $W=G_{10}$ with Coxeter number $h=24$.  We will produce hsops of all degrees less than $h+2=26$ coprime to $h$ (that is, $p=1, 5, 7, 11, 13, 17, 19, 23, 25$).
    Let $\zeta\coloneq e^\frac{2\pi \sqrt{-1}}{12}$. Recall that the reflection representation $V$ of $G_{10}$ is generated by~\cite{JM15,GH96}
    \[
    \renewcommand\arraystretch{1.5}
    s=\left[
\begin{matrix}
1 & 0  \\
0 & \zeta^4
\end{matrix}
\right] \text{ and }
t=\left[
\begin{matrix}
-\frac{2}{3}\zeta^4-\frac{2}{3}\zeta^7-\frac{1}{3}\zeta^8-\frac{1}{3}\zeta^{11} & \frac{1}{6}\zeta^4+\frac{1}{6}\zeta^7-\frac{1}{6}\zeta^8-\frac{1}{6}\zeta^{11}  \\
\frac{1}{3}\zeta^4+\frac{1}{3}\zeta^7-\frac{1}{3}\zeta^8-\frac{1}{3}\zeta^{11} & -\frac{1}{3}\zeta^4-\frac{1}{3}\zeta^7-\frac{2}{3}\zeta^8-\frac{2}{3}\zeta^{11}
\end{matrix}
\right].
    \]
The dual $V^* = \mathrm{span}_\CC\{x,y\}$ gives rise to the polynomial ring $S=\C[x,y]$, and the ring of invariant polynomials $S^W$ is generated by the homogeneous invariant polynomials
\[
f_{12}=x^{12} -704x^3y^9 -88x^9y^3 + -64y^{12} \text{ and }
\]
\[
f_{24}=x^{21}y^3-21x^{18}y^6+123x^{15}y^9-7x^{12}y^{12}-984x^9y^{15}-1344x^6y^{18}-512x^3y^{21},
\]
where the subscript of the invariant polynomials denotes the degree.  Obviously $V^*(-1)$ is carried by $\Theta_1=(x,y)$.

We will illustrate our process by constructing the hsops carrying $V^*(-13)$ and $V^*(-25)$. Beginning with the degree 13 case, we use {\sf Mathematica}~\cite{Mathematica} to construct two arbitrary polynomials, $(x_{13},y_{13})$, with unknown coefficients. For multiple chosen values of $x$ and $y$, we build a system of equations from the condition that the action of $s$ and $t$ on $x_{13}$ and $y_{13}$ is equivalent to the action of $s$ and $t$ on $x$ and $y$ followed by replacing each occurrence of $x$ with $x_{13}$ and $y$ with $y_{13}$. We use this system of equations to solve for our unknown coefficients and compute a degree $13$ hsop $\Theta_{13}=(x_{13},y_{13})$ carrying $V^*(-13)$:
\begin{align*}
x_{13}&=-\frac{1}{14}x^{13} + \frac{109}{7}x^{10}y^3 - \frac{507}{7}x^7y^6 + 28x^4y^9 +64xy^{12},\\
y_{13}&=-x^{12}y + \frac{7}{2}x^9y^4 + \frac{507}{7}x^6y^7 + \frac{872}{7}x^3y^{10} + \frac{32}{7}y^{13}.
\end{align*}
We then use Sage~\cite{sagemath} to verify that the quotient $S/(\Theta_{13})$ was of the correct dimension $13^2 =169$.  A natural candidate for an hsop carrying $V^*(-25)$ is then
\[
\Theta_{25} = (x_{13}f_{12}+xf_{24},y_{13}f_{12}+yf_{24}).
\] We again verified using Sage that $S/(\Theta_{25})$ was of the correct (finite) dimension $25^2= 625$. 

For the remaining hsops of degrees $p=5,7,11,17,19,$ and $23$, we followed the same method, again verifying that the dimension of $S/(\Theta_p)$ was indeed $p^2$ for each hsop listed below. 

\begin{align*}
(x_5,y_5) &= \left(-x^5 - 10x^2y^3 , -5x^3y^2 + 4 y^5\right)\\
  (x_{17},y_{17}) &= \Big(\frac{1}{256}x^{17} - \frac{1}{32}x^{14} y^3 - \frac{465}{64}x^{11} y^6 + \frac{803}{128} x^8 y^9+ \frac{23}{8}x^5 y^{12} + 15 x^2 y^{15}, \\
    &-\frac{15}{128} x^{15} y^2 + \frac{23}{128} x^{12} y^5 - \frac{803}{256}x^9 y^8- \frac{465}{16}x^6 y^{11} + x^3 y^{14} + y^{17}\Big) \\
(x_7,y_7) &= \left(\frac{1}{8}x^7 - \frac{7}{2}x^4y^3 - 7xy^6,-\frac{7}{8} x^6 y + \frac{7}{2}x^3 y^4 + y^7\right)\\
(x_{19},y_{19}) &= \Big(-\frac{1}{512}x^{19} + \frac{1}{64}x^{16} y^3 - \frac{413}{128}x^{13} y^6 - \frac{93}{8} x^{10} y^9- \frac{625}{176}x^7 y^{12} - \frac{1571}{44}x^4 y^{15} - \frac{23}{11}x y^{18},\\
&\frac{23}{5632}x^{18} y - \frac{1571}{2816}x^15 y^4 + \frac{625}{1408}x^{12} y^7 - \frac{93}{8}x^9 y^{10} + \frac{413}{16}x^6 y^{13} + x^3 y^{16} + y^{19}\Big)\\
(x_{11},y_{11})&=\left(-\frac{1}{32}x^{11} + \frac{33}{16}x^8 y^3 + \frac{11}{2}x^2 y^9, \frac{11}{32}x^9 y^2 + \frac{33}{4}x^3 y^8 + y^{11}\right)\\
(x_{23},y_{23}) &= \Big(\frac{1}{2048}x^{23} - \frac{1}{256}x^{20} y^3 + \frac{795}{512}x^{17} y^6+ \frac{41815}{7168}x^{14} y^9 \\ &+ \frac{7663}{256}x^{11} y^{12} - \frac{25089}{896}x^8 y^{15} + \frac{265}{8}x^5y^{18} + \frac{2}{7}x^2 y^{21},\\
&-\frac{1}{3584}x^{21} y^2 + \frac{265}{1024}x^{18} y^5 + \frac{25089}{14336}x^{15} y^8 \\ &+ \frac{7663}{512}x^{12} y^{11} - \frac{41815}{1792}x^9 y^{14} + \frac{795}{16}x^6 y^{17} + 
 x^3 y^{20} + y^{23}\Big).
\end{align*}

\end{example}

We extend the definition of the algebraic parking space to all irreducible well--generated complex reflection groups. Let $x_1,\ldots,x_n$ be a basis for $V^*$ and let $\Theta_{h+1}$ be an hsop of degree $h+1$ carrying $V^*(-h-1)$ such that the map $x_i\mapsto \theta_i$ is $W$-equivariant.  
\begin{definition}
    Let $W$ be an irreducible well--generated complex reflection group. Let $\Theta_{h+1}=(\theta_1,\ldots,\theta_n)$ and ${\bf x}=(x_1,\ldots,x_n)$ be as above. Define the \defn{algebraic $W$-parking space} by 
    \[
    \Park^{alg}(W,\Theta_{h+1})\coloneq S/(\Theta_{h+1}-\textbf{x}).
    \]
    \end{definition}
    We similarly define the algebraic parking space for the Fuss case. 
    \begin{definition}
    Let $W$ be an irreducible well--generated complex reflection group and let $k\geq 1$. Let $\Theta_{kh+1}=(\theta_1,\ldots,\theta_n)$ and ${\bf x}=(x_1,\ldots,x_n)$ be as above. Define the \defn{$k$-algebraic $W$-parking space} by 
    \[
    \Park_{alg}^k(W,\Theta_{kh+1})\coloneq S/(\Theta_{kh+1} - {\bf x}).
    \]
\end{definition}
The ideals $(\Theta_{h+1}-\textbf{x})$ and $(\Theta_{kh+1}-\textbf{x})$ are $(W\times C)$-stable and $(W\times \Z_{kh})$-stable, respectively, with the actions of linear substitutions and scalar substitutions identical to the real case. We may drop the dependence on the hsop chosen by extending the proofs of~\Cref{prop:alg_park_iso_btw_hsops} and~\Cref{prop:k_alg_park_iso_btw_hsops}. 

\subsection{Results}
Our main results exhibit a $(W\times C)$-equivariant (or $(W\times \Z_{kh})$-equivariant in the Fuss case) isomorphism between the two parking spaces and enumerate the noncrossing parking functions for all irreducible well--generated complex reflection groups. Let $W$ be an irreducible well--generated complex reflection group not of type $G_{34}$, $E_7$, or $E_8$. Let $h$ be the Coxeter number of $W$ and let $n$ be its rank. 

\begin{theorem}
    \label{thm: bijection}
    For all such $W$, there is a $(W\times C)$-equivariant isomorphism between $\Park^\NC(W)$ and $\Park^{alg}(W)$.
\end{theorem}

Combining~\Cref{thm: bijection} with computer calculations using~\cite{JM15,GH96}, we obtain the following enumeration:
\begin{theorem}
    \label{thm: Park is (h+1)^n}
    For all such $W$, there are $(h+1)^n$ $W$-noncrossing parking functions. 
\end{theorem}

We will see~\Cref{thm: bijection} implies~\Cref{thm: Park is (h+1)^n} for $G(d,1,n)$ and $G(d,d,n)$ by the construction of our algebraic parking space. We do the same in the Fuss case. Let $k\geq 1$. 

\begin{theorem}
\label{thm: k-bijection}
    For all such $W$, there is a $(W\times \Z_{kh})$-equivariant isomorphism between $\Park^k_\NC(W)$ and $\Park^k_{alg}(W)$.
\end{theorem}
Again using~\cite{JM15,GH96} and~\Cref{thm: k-bijection}, we obtain the following:
\begin{theorem}
\label{thm: k-Park is (h+1)^n}
    For all such $W$, there are $(kh+1)^n$ $k$-$W$-noncrossing parking functions. 
\end{theorem}

Again, we will see that~\Cref{thm: k-bijection} implies~\Cref{thm: k-Park is (h+1)^n} by our choice of hsop in each type. We prove~\Cref{thm: bijection} separately for $G(d,1,n)$ and $G(d,d,n)$ in section~\Cref{sec: d1n bijection} and~\Cref{sec:ddn bijection}, respectively. Similarly, we prove~\Cref{thm: k-bijection} for $G(d,1,n)$ in~\Cref{sec:fuss_d1n_bij} and $G(d,d,n)$ in section~\Cref{sec:fuss_ddn_bij}. For the exceptional groups, we enumerate the noncrossing parking functions using~\cite{JM15,GH96}. We then compute the character of $W\times C$ acting on $S/(\Theta_p)$ using~\cite{JM15,GH96} and extend Proposition 2.15 in~\cite{armstrong2015parking} to our case to verify the isomorphism. We do the same in the Fuss case extending Proposition 3.5 and Equation 1.7 in~\cite{rhoades2014parking}. 

\section{\texorpdfstring{$G(d,1,n)$}{Type G(d,1,n)}}
\label{sec:gd1n}
The irreducible well--generated complex reflection group $G(d,1,n)$ is the group of $n\times n$ monomial matrices with nonzero entries $h_i\in\mu_d$. We may think of the elements of this group as permutations of the set
\[
[n]^d\coloneq\{1^1,2^1,\ldots,n^1,1^2,\ldots,n^2,\ldots,1^d,\ldots,n^d\}.
\]
For any $i^j\in [n]^d$, we will call $i$ its \defn{number} and $j$ its \defn{color}. We first recall a combinatorial model for the $G(d,1,n)$-noncrossing partitions~\cite{bessis2006eer}. We will then use these to construct the noncrossing parking functions in this type. Since $G(2,1,n)$ recovers the Weyl group of type $B_n$, we are able to extend isomorphisms in type $B$ from~\cite{armstrong2015parking} to our case. Take the Coxeter element $c = [1,2,3,\ldots,n]$, where the square brackets indicate $i^j$ maps to $(i+1)^j$ for all $1\leq i\leq n-1$ and $n^j$ maps to $1^{j+1\mod d}$. We use double parenthesis in cycle notation to indicate a cycle occurring in all $d$ colors.

\subsection{Combinatorial Model for \texorpdfstring{$G(d,1,n)$}{G(d,1,n)}}
\label{sec: d1n model}
In~\cite{bessis2006eer}, it is shown that the $G(d,1,n)$-noncrossing partitions are the type $A$ noncrossing partitions that are invariant under $\frac{2\pi}{n}$ rotation. As such, we can visualize the noncrossing partitions in $G(d,1,n)$ by taking $dn$ vertices around a circle and identifying them with the set
\[
[n]^d\coloneq\{1^1,\ldots,n^1,1^2,\ldots,n^2,\ldots,1^d,\ldots,n^d\}.
\]
We label the $dn$ vertices by the elements of $[n]^d$ clockwise increasing in number and clockwise decreasing in color. We identify a set partition of the set $[n]^d$, $\pi=\{B_1,\dots,B_r\}$, with the convex hulls of each block.  Given such a partition $\pi$, the \defn{zero block} (if it exists) is the unique block, denoted $B_0$, satisfying the condition that if $i^r \in B_0$ for any $r=1,\dots,d$, then $i^j \in B_0$ for all $j=1,\dots,d$. 

\begin{example}
In~\Cref{fig:noncrossing partition construction} we see an example of a noncrossing partition with zero block $\{6^1,6^2,6^3\}$ given by the partition 
\begin{equation}
\label{partition pi for d,1,n example}
    \pi=\{\{1^1,4^1,5^1\},\{2^1,3^1\},\{1^2,4^2,5^2\},\{2^2,3^2\},\{1^3,4^3,5^3\},\{2^3,3^3\},\{6^1,6^2,6^3\}\}.
\end{equation}  
\begin{figure}[h]
    \centering
    \begin{tikzpicture}
%Connect zero block 
\draw[fill=lightgray] (-360/18*9:-3cm)--
(-360/18*3:-3cm) --
(-360/18*15:-3cm)--
(-360/18*9:-3cm);

%Connect the 2-3
\draw[] (-360/18*5:-3cm)--
(-360/18*6:-3cm);

\draw[] (-360/18*11:-3cm)--
(-360/18*12:-3cm);

\draw[] (-360/18*17:-3cm)--
(-360/18*18:-3cm);

%Connect the 1-4-5
\draw[fill=lightgray] (-360/18*4:-3cm)--
(-360/18*7:-3cm) --
(-360/18*8:-3cm)--
(-360/18*4:-3cm);

\draw[fill=lightgray] (-360/18*10:-3cm)--
(-360/18*13:-3cm) --
(-360/18*14:-3cm)--
(-360/18*10:-3cm);

\draw[fill=lightgray] (-360/18*16:-3cm)--
(-360/18*1:-3cm) --
(-360/18*2:-3cm)--
(-360/18*16:-3cm);

\foreach \a/\n in {1^1/4,2^1/5,3^1/6,4^1/7,5^1/8,6^1/9,1^3/10,2^3/11,3^3/12,4^3/13,5^3/14,6^3/15,1^2/16,2^2/17,3^2/18,4^2/1,5^2/2,6^2/3}
{
\draw[fill](-360/18*\n:-3cm)circle(1.2pt)
node[anchor=-\n*360/18] (\a) {$\a$};
}

\end{tikzpicture}
    \caption{A $G(3,1,6)$-noncrossing partition with $\pi$ as in~\ref{partition pi for d,1,n example}.
    }
    \label{fig:noncrossing partition construction}
\end{figure}
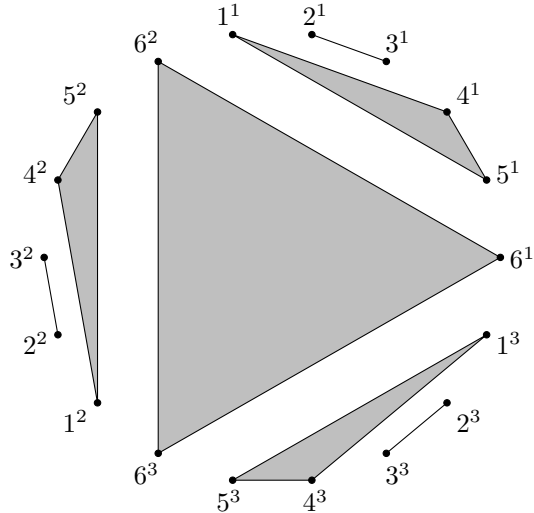
\end{example}

We can visualize a $G(d,1,n)$-noncrossing parking function, $[w,\pi]$, by extending the labelling in types $A$, $B$, and $D$ from~\Cref{subsec: ncp_in_rrg} to this context; taking noncrossing partition $\pi$ and labelling each block $B_i \in \pi$ by $w(B_i)$.

\begin{example}\label{ex:pf_g316}
Continuing with the partition $\pi$ as in ~\Cref{fig:noncrossing partition construction}, we have an element $[w,\pi]\in \Park^\NC(W)$ given by
\begin{equation}
\label{eqn:parking-function-example}
\begin{aligned}
wW_\pi&=
\left(
\begin{matrix}
1^1 & 2^1 & 3^1 & 4^1 & 5^1 & 6^1 \\
2^2 & 3^3 & 1^1 & 5^2 & 6^2 & 4^2
\end{matrix}
\right)W_\pi \\
&=
\left(
\begin{array}{ccc|cc|c}
1^1 & 4^1 & 5^1 & 2^1 & 3^1 & 6^1\\
2^2 & 5^2 & 6^2 & 3^3 & 1^1 & 4^2 
\end{array}
\right)W_\pi.
\end{aligned}
\end{equation}
We model this noncrossing parking function in~\Cref{fig:pf_example_gd1n}.
\end{example}
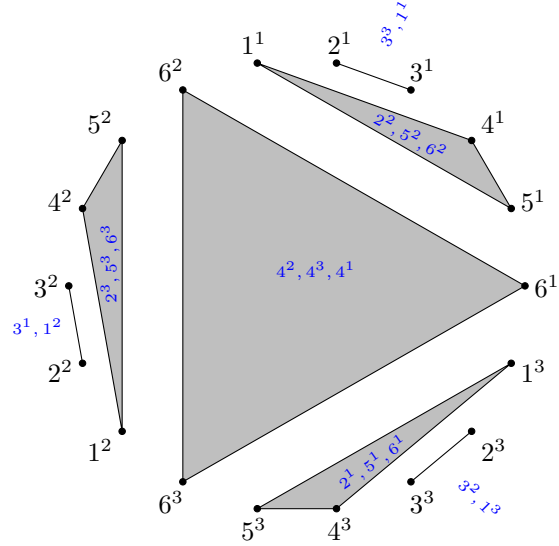
\begin{figure}[h]
    \centering
        \begin{tikzpicture}
%Connect zero block 
\draw[fill=lightgray] (-360/18*9:-3cm)--node[below,yshift=-.85cm,xshift=-.5cm]{\tiny \color{blue}$4^2,4^3,4^1$}
(-360/18*3:-3cm) --
(-360/18*15:-3cm)--
(-360/18*9:-3cm);

%Connect the 2-3
\draw[] (-360/18*5:-3cm)--node[above,rotate=60,xshift=.75cm,yshift=-.15cm]{\tiny \color{blue}$3^3,1^1$}
(-360/18*6:-3cm);

\draw[] (-360/18*11:-3cm)--node[below,rotate=-45,xshift=.75cm,yshift=.2cm]{\tiny \color{blue} $3^2,1^3$}
(-360/18*12:-3cm);

\draw[] (-360/18*17:-3cm)--node[below,xshift=-.5cm,yshift=.2cm]{\tiny \color{blue}$3^1,1^2$}
(-360/18*18:-3cm);

%Connect the 1-4-5
\draw[fill=lightgray] (-360/18*4:-3cm)--node[below,rotate=-30,yshift=.15cm,xshift=.75cm]{\tiny \color{blue}$2^2,5^2,6^2$}
(-360/18*7:-3cm) --
(-360/18*8:-3cm)--
(-360/18*4:-3cm);

\draw[fill=lightgray] (-360/18*10:-3cm)--node[above, rotate=35,xshift=-.75cm,yshift=-.15cm]{\tiny \color{blue} $2^1,5^1,6^1$}
(-360/18*13:-3cm) --
(-360/18*14:-3cm)--
(-360/18*10:-3cm);

\draw[fill=lightgray] (-360/18*16:-3cm)--node[above,rotate=90,yshift=-.35cm,xshift=.75cm]{\tiny \color{blue} $2^3,5^3,6^3$}
(-360/18*1:-3cm) --
(-360/18*2:-3cm)--
(-360/18*16:-3cm);

\foreach \a/\n in {1^1/4,2^1/5,3^1/6,4^1/7,5^1/8,6^1/9,1^3/10,2^3/11,3^3/12,4^3/13,5^3/14,6^3/15,1^2/16,2^2/17,3^2/18,4^2/1,5^2/2,6^2/3}
{
\draw[fill](-360/18*\n:-3cm)circle(1.2pt)
node[anchor=-\n*360/18] (\a) {$\a$};
}

\end{tikzpicture}
    \caption{A $G(3,1,6)$-parking function example with partition $\pi$ as in~\ref{partition pi for d,1,n example} and $w$ from~\ref{eqn:parking-function-example}.}
    \label{fig:pf_example_gd1n}
\end{figure}

\begin{figure}
    \centering
        \begin{tikzpicture}[x=1.0cm,y=1.0cm,scale=1,baseline={(0,0)}]
%Connect zero block 
\draw[fill=lightgray] (-360/18*9:-3cm)--node[below,yshift=-.85cm,xshift=-.5cm]{\tiny \color{blue}$1^2,1^3,1^1$}
(-360/18*3:-3cm) --
(-360/18*15:-3cm)--
(-360/18*9:-3cm);

%Connect the 2-3
\draw[] (-360/18*5:-3cm)--node[above,rotate=60,xshift=.75cm,yshift=-.15cm]{\tiny \color{blue}$2^2,4^1$}
(-360/18*6:-3cm);

\draw[] (-360/18*11:-3cm)--node[below,rotate=-45,xshift=.75cm,yshift=.2cm]{\tiny \color{blue} $2^1,4^3$}
(-360/18*12:-3cm);

\draw[] (-360/18*17:-3cm)--node[below,xshift=-.5cm,yshift=.2cm]{\tiny \color{blue}$2^3,4^2$}
(-360/18*18:-3cm);

%Connect the 1-4-5
\draw[fill=lightgray] (-360/18*4:-3cm)--node[below,rotate=-30,yshift=.15cm,xshift=.75cm]{\tiny \color{blue}$3^3,5^2,6^2$}
(-360/18*7:-3cm) --
(-360/18*8:-3cm)--
(-360/18*4:-3cm);

\draw[fill=lightgray] (-360/18*10:-3cm)--node[above, rotate=35,xshift=-.75cm,yshift=-.15cm]{\tiny \color{blue} $3^2,5^1,6^1$}
(-360/18*13:-3cm) --
(-360/18*14:-3cm)--
(-360/18*10:-3cm);

\draw[fill=lightgray] (-360/18*16:-3cm)--node[above,rotate=90,yshift=-.35cm,xshift=.75cm]{\tiny \color{blue} $3^1,5^3,6^3$}
(-360/18*1:-3cm) --
(-360/18*2:-3cm)--
(-360/18*16:-3cm);

\foreach \a/\n in {1^1/4,2^1/5,3^1/6,4^1/7,5^1/8,6^1/9,1^3/10,2^3/11,3^3/12,4^3/13,5^3/14,6^3/15,1^2/16,2^2/17,3^2/18,4^2/1,5^2/2,6^2/3}
{
\draw[fill](-360/18*\n:-3cm)circle(1.2pt)
node[anchor=-\n*360/18] (\a) {$\a$};
}

\end{tikzpicture}
\qquad
        \begin{tikzpicture}[x=1.0cm,y=1.0cm,scale=1,baseline={(0,0)}]
%Connect zero block 
\draw[fill=lightgray] (-360/18*10:-3cm)--node[below,yshift=-.8cm,xshift=-.9cm]{\tiny \color{blue}$4^2,4^3,4^1$}
(-360/18*4:-3cm) --
(-360/18*16:-3cm)--
(-360/18*10:-3cm);

%Connect the 2-3
\draw[] (-360/18*6:-3cm)--node[above,rotate=50,xshift=.75cm,yshift=-.25cm]{\tiny \color{blue}$3^3,1^1$}
(-360/18*7:-3cm);

\draw[] (-360/18*12:-3cm)--node[below,rotate=-65,xshift=.75cm,yshift=.2cm]{\tiny \color{blue} $3^2,1^3$}
(-360/18*13:-3cm);

\draw[] (-360/18*18:-3cm)--node[below,rotate=-10,xshift=-.5cm,yshift=.2cm]{\tiny \color{blue}$3^1,1^2$}
(-360/18*1:-3cm);

%Connect the 1-4-5
\draw[fill=lightgray] (-360/18*5:-3cm)--node[below,rotate=-50,yshift=.15cm,xshift=.75cm]{\tiny \color{blue}$2^2,5^2,6^2$}
(-360/18*8:-3cm) --
(-360/18*9:-3cm)--
(-360/18*5:-3cm);

\draw[fill=lightgray] (-360/18*11:-3cm)--node[above,rotate=15,xshift=-.75cm,yshift=-.15cm]{\tiny \color{blue} $2^1,5^1,6^1$}
(-360/18*14:-3cm) --
(-360/18*15:-3cm)--
(-360/18*11:-3cm);

\draw[fill=lightgray] (-360/18*17:-3cm)--node[above,rotate=70,yshift=-.35cm,xshift=.75cm]{\tiny \color{blue} $2^3,5^3,6^3$}
(-360/18*2:-3cm) --
(-360/18*3:-3cm)--
(-360/18*17:-3cm);

\foreach \a/\n in {1^1/4,2^1/5,3^1/6,4^1/7,5^1/8,6^1/9,1^3/10,2^3/11,3^3/12,4^3/13,5^3/14,6^3/15,1^2/16,2^2/17,3^2/18,4^2/1,5^2/2,6^2/3}
{
\draw[fill](-360/18*\n:-3cm)circle(1.2pt)
node[anchor=-\n*360/18] (\a) {$\a$};
}

\end{tikzpicture}
    \caption{The image of the $G(3,1,6)$-noncrossing parking function in~\Cref{ex:pf_g316} under the action of $v\in W$ as in~\Cref{ex:gd1n_wxc_action} (left) and $c\in C$ (right).}
    \label{fig:gd1n_pf_wxc_action}
\end{figure}
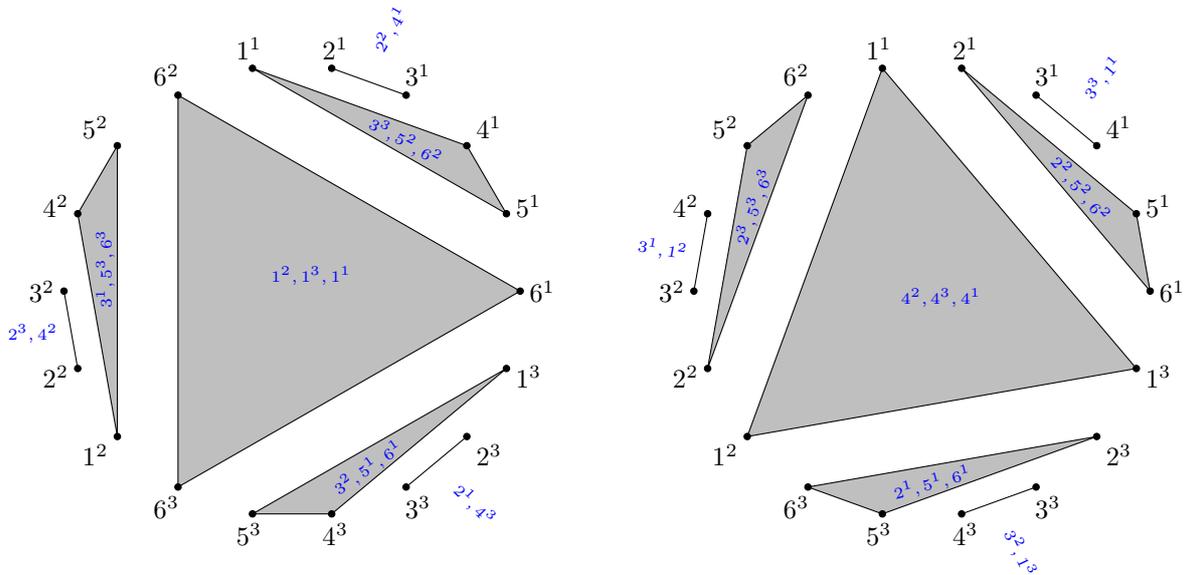

The $(W\times C)$-action in this setting corresponds to $W$ permuting the labels on each block of $\pi$ and $C$ rotating the picture $\frac{2\pi}{dn}$ clockwise. 
\begin{example}\label{ex:gd1n_wxc_action}
    Consider the $G(3,1,6)$-parking function in~\Cref{fig:pf_example_gd1n}. We have Coxeter element $c=[1,2,3,4,5,6]$. Define $\Omega(\pi)$ for a noncrossing partition $\pi$ of $[n]^d$ by sending $\pi$ to its the corresponding word with cycles given by the blocks of $\pi$. We abuse notation and sometimes write $c\pi c^{-1}$ to denote $c\Omega(\pi) c^{-1}$. Consider $v=(\!(1 4)\!)(\!(2 \overline{3})\!)(\!(5 6)\!)\in G(3,1,6)$ where the double parenthesis indicate the cycle occurring in all $d$ colors and $i=i^1$, $\overline{i}=i^2$, and $\overline{\overline{i}}=i^3$ for all $i\in[6]$. The models in~\Cref{fig:gd1n_pf_wxc_action} show the action of $v$ (left) and the action of $c\in C$ (right). The figure on the left shows the noncrossing parking function $[vw,\pi]$ and the figure on the right shows $[wc^{-1},c\pi c^{-1}]$.
\end{example}

\subsection{Proof of~\Cref{thm: bijection} for \texorpdfstring{$G(d,1,n)$}{G(d,1,n)}}\label{sec: d1n bijection}
In this section we extend the isomorphism in type $B$ from~\cite{armstrong2015parking} to $G(d,1,n)$. Let $W=G(d,1,n)$ act on $V=\C^n$ with the standard coordinate functionals ${\bf x}=x_1,\dots,x_n$ and choose 
\[
(\Theta_{h+1}) = (\theta_1,\dots,\theta_n)\coloneq (x_1^{dn+1},\dots,x_n^{dn+1}).
\]
By~\Cref{rmk:hsop_inf_fams}, $\Theta_{h+1}$ is an hsop of degree $dn+1=h+1$ carrying $V^*(-h-1)$, where $h=dn$ is the Coxeter number of $W$. Further note that the map sending $x_i \to \theta_i$ is $W$-equivariant. Now consider the ideal
\begin{equation*}
\begin{split}
    (\Theta_{h+1}-\textbf{x}) &= (x_1^{dn+1}-x_1,\dots,x_n^{dn+1}-x_n) \\
&=(x_1(x_1^{dn}-1),\dots,x_n(x_n^{dn}-1)).
\end{split}
\end{equation*}
The subvariety cut out by $(\Theta_{h+1} - {\bf x})$ is
\[
V^{\Theta_{h+1}} \coloneq\{(v_1,\dots,v_n)\in\C^n : v_i=0 \text{ or } v_i^{dn}=1\}.
\] 
Specifically, each $v_i\in V^{\Theta_{h+1}}$ is of the form $\lambda^j\omega^i$ for some $i,j \in \mathbb{N}$, where $\omega\coloneq e^{\frac{2\pi \sqrt{-1}}{n}},$ and $\lambda\coloneq e^{\frac{2\pi \sqrt{-1}}{d}}$. Note also $|V^{\Theta_{h+1}}|=(h+1)^n$. We first describe the forward isomorphism
\[
\Park^{\NC}(W) \xrightarrow{f} V^{\Theta_{h+1}}.
\]

Begin with a noncrossing partition $\pi$ of $[n]^d$ as described in ~\Cref{sec: d1n model}. We then re-encode the given partition via a bijection from~\cite{reiner1997classical} by a parenthesization of the doubly infinite string
\[
\dots 1^1,2^1,\dots,(n-1)^1,n^1,1^d,\dots,(n-1)^d, n^d, 1^{d-1},2^{d-1},\dots,n^{d-1},\dots,1^2,2^2,\dots,n^2,1^1,2^1\dots.
\]

The idea behind this re-encoding is that each nonzero block $B_i$ of $\pi$ corresponds to a pair $(L,R)$, with $L,R\subset [n]^d$, of left parenthesis and right parenthesis. If $r\in L$, place a left parenthesis to the left of any $r^j$ in the doubly infinite string for all $j=1,\ldots,d$. Similarly, place right parenthesis using elements in $R$. Each pair of left and right parenthesis corresponds to a block $B$ containing the same elements as the pair of parentheses does, excluding any elements contained in a nested pair of parentheses. Any elements not belonging to a parentheses pair belong to the zero block. 
\begin{example}\label{ex:gd1n}
The noncrossing parking function in~\Cref{fig:pf_example_gd1n} produces the parenthesization
\[
\dots,(1^1,(2^1,3^1),4^1,5^1),6^1,(1^3,(2^3,3^3),4^3,5^3),6^3,
(1^2,(2^2,3^2),4^2,5^2),6^2,(1^1,\dots
\]
corresponding to the pair $(L,R) = (\{1,2\},\{3,5\})$. The unparenthesized $6^i$'s correspond to the zero block $B_0 = \{6^1,6^2,6^3\}$.
\end{example}
\begin{definition}
    In a given pair of parentheses, let $i^j$ be the element appearing immediately after the left parenthesis. We will call $i^j$ the \defn{opener } of the corresponding block $B$.
\end{definition}

We define our isomorphism by the map $f$ which sends $[w,\pi]$ to the vector $v = (v_1,\dots,v_n)\in V^{\Theta_{h+1}}$ with
\begin{itemize}
    \item $v_r = 0$ if $r^1\in w(B_0)$ with $B_0$ the zero block, or
    \item $v_r = \lambda^j\omega^i$ if $r^1\in w(B)$ with $B$ having opener $i^j$,
\end{itemize}
where $\lambda \coloneq e^\frac{2\pi \sqrt{-1}}{d}$ and $\omega \coloneq e^\frac{2\pi \sqrt{-1}}{n}$.
\begin{example}
We continue with~\Cref{ex:gd1n}. Here, the zero block $B_0 = \{6^1,6^2,6^3\}$ has image $w(B_0) = \{4^1,4^2,4^3\}$, so $v_4 = 0$. We tabulate the nonzero blocks, their openers, and images:
\begin{center}
\bgroup
\def\arraystretch{1.5}
\begin{tabular}{|c|c|c|}\hline
$B_i$ & opener & $w(B_i)$ \\\hline\hline
$\{2^1,3^1\}$     & $2^1$ & $\{3^3,1^1\}$     \\\hline
$\{1^1,4^1,5^1\}$ & $1^1$ & $\{2^2,5^2,6^2\}$ \\\hline
$\{2^2,3^2\}$     & $2^2$ & $\{3^1,1^2\}$  \\\hline
$\{1^2,4^2,5^2\}$ & $1^2$ & $\{2^3,5^3,6^3\}$    \\\hline
$\{2^3,3^3\}$     & $2^3$ & $\{3^2,1^3\}$ \\\hline
$\{1^3,4^3,5^3\}$ & $1^3$ & $\{2^1,5^1,6^1\}$ \\\hline
\end{tabular}.
\egroup
\end{center}

So $1^1 \in \{3^3,1^1\}$ with opener $2^1$, and thus $v_1 = \lambda^1\omega^2$. Continuing in this way, we see
\begin{equation}
\label{eqn:vector image of parking example d1n}
    v=(\lambda^1\omega^2,\lambda^3\omega^1,\lambda^2\omega^2,0,\lambda^3\omega^1,\lambda^3\omega^1).
\end{equation}
\end{example}

To show this map is indeed a bijection, we now describe the inverse \[
V^{\Theta_{h+1}} \xrightarrow{f^{-1}} \Park^{\NC}(W).
\] 

Given a vector $v \in V^{\Theta_{h+1}}$, we can determine the parenthesization of the infinite string by noting that there will be a left parenthesis located immediately preceding $r^1,r^2,\dots,r^d$ if and only if at least one coordinate of $v$  is equal to $\lambda^j\omega^r$ for some $j$. Next, we can recover the right parenthesis by finding the number of elements contained in the parentheses pair. Let $m_r$ be the number of occurrences of $\lambda^j\omega^r$ for any $j=1,\dots,d$, amongst the coordinates of $v$. Then, given any such sequence of multiplicities $(m_1,\dots,m_n)$, beginning with the smallest $m_r$ to the largest, we recover the right parenthesis by placing it immediately after the unique element in the infinite string so that the given pair of parenthesis contains $m_r$ elements, excluding all elements previously enclosed in a pair of parentheses. The zero block again corresponds to the unparenthesized letters in the infinite string. 

\begin{example}
\label{ex:d1n_inverse_partition}
Consider the vector from~\ref{eqn:vector image of parking example d1n}. We first see that $\lambda^1\omega^2$ and $\lambda^2\omega^2$ occur, giving a left parenthesis in the infinite string immediately before each $2^j$ for $j=1,2,3$. Similarly $\lambda^3\omega^1$ occurs three times, giving a left parenthesis immediately before each $1^j$ for $j=1,2,3$. So we place our left parentheses as such:
\[
\dots,(1^1,(2^1,3^1,4^1,5^1,6^1,(1^3,(2^3,3^3,4^3,5^3,6^3,
(1^2,(2^2,3^2,4^2,5^2,6^2,(1^1,\dots.
\]

Next, we list our multiplicities:
\begin{equation*}
\begin{split}
    m_1 &= 3 \\
    m_2 &= 2 \\
    m_3 &= m_4=m_5=m_6=0.
\end{split}
\end{equation*}

Given this, there is a right parenthesis occurring after each 2 elements following a left parenthesis that has opener $2^j$ for $j=1,2,3$. So we place the right parenthesis as such:

\[
\dots,(1^1,(2^1,3^1),4^1,5^1,6^1,(1^3,(2^3,3^3),4^3,5^3,6^3,
(1^2,(2^2,3^2),4^2,5^2,6^2,(1^1,\dots.
\]

Next, ignoring $2^j,3^j$ for $j=1,2,3$, there is a right parenthesis occuring after each 3 elements following a left parenthesis that has opener $1^j$ for $j=1,2,3$. Now our infinite string becomes 

\[
\dots,(1^1,(2^1,3^1),4^1,5^1),6^1,(1^3,(2^3,3^3),4^3,5^3),6^3,
(1^2,(2^2,3^2),4^2,5^2),6^2,(1^1,\dots.
\]
\end{example}

After obtaining the infinite parenthesization, we have the partition $\pi$ of $[n]^d$. We have only left to determine the labellings of each block $B$. Suppose $B$ is opened by $i^j$. Then, $r^1 \in [n]^d$ lies in the set $w(B)$ if and only if $v_r = \lambda^j\omega^i$. Then, as $w\in G(d,1,n)$, we recover the labellings for each $r^j\in [n]^d$. 

\begin{example}
\label{ex:d1n_inverse_labelling}
Continuing with~\Cref{ex:d1n_inverse_partition}, $v_1=\lambda^1\omega^2$, giving that $1^1\in w(\{2^1,3^1\})$. From here, we see that as $w$ maps $2^1$ or $3^1$ to $1^1$. Thus, $w$ maps $2^2$ or $3^2$ to $1^2$, and so $1^2\in w(\{2^2,3^2\})$. In the same way, $1^3\in w(\{2^3,3^3\})$. For the remainder of the labellings, we have
\begin{itemize}
    \item $2^1 \in w(\{1^3,4^3,5^3\})$,
    \item $3^1\in w(\{2^2,3^2\})$,
    \item $5^1 \in w(\{1^3,4^3,5^3\})$,
    \item $6^1 \in w(\{1^3,4^3,5^3\})$.
\end{itemize}
Analogously to $1^j$, we recover the labellings for $2^j$, $3^j$, $5^j$, and $6^j$ with $j=2,3$.
\end{example}

Then we have recovered $[w,\pi] = f^{-1}(v)$ and shown this map is a bijection. We leave it to the reader to verify this bijection is $(W\times C)$-equivariant. This gives the proof of~\Cref{thm: bijection} in this case and~\Cref{thm: Park is (h+1)^n} by our choice of hsop.

\section{Fuss Analogues of Type \texorpdfstring{$G(d,1,n)$}{G(d,1,n)}}
\label{sec:fuss_gd1n}
In~\cite{rhoades2014parking}, B.~Rhoades goes through the construction of a combinatorial model for a Fuss noncrossing parking function in type $B$.  We can adapt the same maps to our models in $W=G(d,1,n)$ to visualize the $k$-$W$-noncrossing partitions.

\subsection{Visualizing Fuss \texorpdfstring{$G(d,1,n)$}{G(d,1,n)}-noncrossing parking functions}
\label{sec:fussd1n_visualize}
We can naively visualize a $k$-$W$-noncrossing parking function by taking $k$ noncrossing partition models from the multichain $(\pi_1\leq\cdots\leq \pi_k)$, decorating $\pi_1$ with images under $w\in W$. Note the labelling on $\pi_1$ induces the labelling on $\pi_i$ for all $i>1$.

\begin{example}
Consider the multichain 
\[
(\{\{1^1,2^1\},\{1^2,2^2\},\{1^3,2^3\}\}\leq\{\{1^1,2^1,3^1\},\{1^2,2^2,3^2\},\{1^3,2^3,3^3\}\})
\]
        and \begin{equation*}
\begin{aligned}
w&=
\left(
\begin{matrix}
1^1 & 2^1 & 3^1 \\
2^1 & 1^2 & 3^2 
\end{matrix}
\right).
\end{aligned}
\end{equation*}
Our naive model is shown in~\Cref{fig:rudimentary fuss model}.
\begin{figure}[h]
        \begin{center}
        \[
        \left(
\tikzset{every node/.style={font=\scriptsize}}
\begin{tikzpicture}[baseline=-1mm]
\draw[](-360/9*1:-.75cm)--node[above,rotate=-75,yshift=-.35cm,xshift=-.6cm]{\tiny\color{blue}$1^2,2^1$}
(-360/9*2:-.75cm);

\draw[](-360/9*4:-.75cm)--node[above,xshift=.6cm,yshift=-.35cm]{\tiny\color{blue}$1^1,2^3$}
(-360/9*5:-.75cm);

\draw[](-360/9*7:-.75cm)--node[below,rotate=60,yshift=.35cm,xshift=-.5cm]{\tiny\color{blue}$1^3,2^2$}
(-360/9*8:-.75cm);

\draw[](-360/9*3:-.75cm)--node[below]{\tiny\color{blue}$3^2$}
(-360/9*3:-.75cm);

\draw[](-360/9*6:-.75cm)--node[above]{\tiny\color{blue}$3^1$}
(-360/9*6:-.75cm);

\draw[](-360/9*9:-.75cm)--node[above,xshift=.3cm,yshift=-.25cm]{\tiny\color{blue}$3^3$}
(-360/9*9:-.75cm);

\foreach \a/\n in {1^1/1,2^1/2,3^1/3,1^3/4,2^3/5,3^3/6,1^2/7,2^2/8,3^2/9}
{
\draw[fill](-360/9*\n:-.75cm)circle(1pt)
node[anchor=-\n*360/9] (\a) {$\a$};
}

\end{tikzpicture},
\qquad
\tikzset{every node/.style={font=\scriptsize}}
\begin{tikzpicture}[baseline=-1mm]

\draw[fill=lightgray](-360/9*1:-.75cm)--
(-360/9*2:-.75cm)--
(-360/9*3:-.75cm)--
(-360/9*1:-.75cm);

\draw[fill=lightgray](-360/9*4:-.75cm)--
(-360/9*5:-.75cm)--
(-360/9*6:-.75cm)--
(-360/9*4:-.75cm);

\draw[fill=lightgray](-360/9*7:-.75cm)--
(-360/9*8:-.75cm)--
(-360/9*9:-.75cm)--
(-360/9*7:-.75cm);

\foreach \a/\n in {1^1/1,2^1/2,3^1/3,1^3/4,2^3/5,3^3/6,1^2/7,2^2/8,3^2/9}
{
\draw[fill](-360/9*\n:-.75cm)circle(1pt)
node[anchor=-\n*360/9] (\a) {$\a$};
}

\end{tikzpicture}
\right)
\]
\end{center}
    \caption{Initial Fuss model for $G(3,1,3)$ with $k=2$.}
    \label{fig:rudimentary fuss model}
\end{figure}
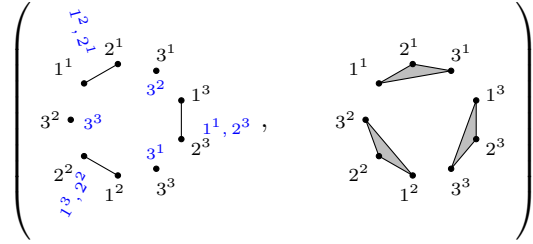
\end{example}

We would like to have one model representing a $k$-$W$-noncrossing parking function, rather $k$ models as in~\Cref{fig:rudimentary fuss model}. To do so we follow the same construction as B.~Rhoades in type B~\cite{rhoades2014parking}.

First, define $\Omega$ which takes a set partition $\pi$ of $[n]^d$ and outputs $\Omega(\pi)\in G(d,1,n)$ by sending the blocks in $\pi$ to cycles with elements in the same order as the infinite parenthesization, as introduced in~\Cref{ex:gd1n_wxc_action}. 
\begin{example}
For example, 
\[
\Omega(\{\{1^1,2^1,3^1\},\{4^1,6^1\}\{5^1\},\{1^2,2^2,3^2\},\{4^2,6^2\}\{5^2\},\{1^3,2^3,3^3\},\{4^3,6^3\}\{5^3\}\}) = (\!(123)\!)(\!(46)\!),
\]
with notation as in~\Cref{ex:gd1n_wxc_action}. 
\end{example}
Next, we define the map $\delta$ on a multichain by
\[
\delta(w_1,w_2,\ldots,w_k) = (w_1^{-1}w_2,w_2^{-1}w_3,\ldots,w_k^{-1}c).
\]

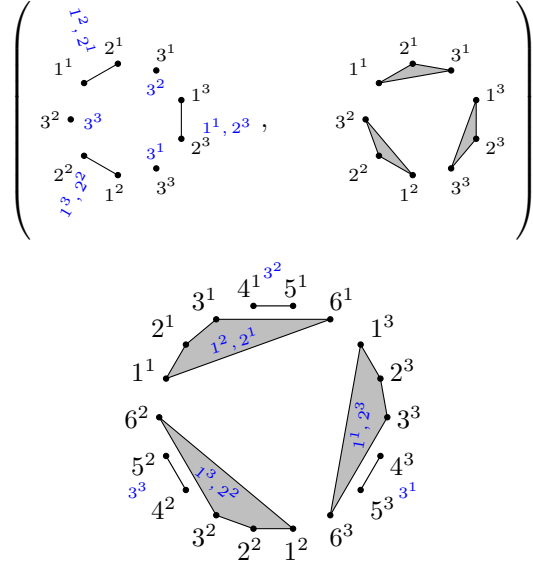
\begin{figure}[h]
    \centering
        \[
        \left(
\tikzset{every node/.style={font=\scriptsize}}
\begin{tikzpicture}[baseline=-1mm]
\draw[](-360/9*1:-.75cm)--node[above,rotate=-75,yshift=-.35cm,xshift=-.6cm]{\tiny\color{blue}$1^2,2^1$}
(-360/9*2:-.75cm);

\draw[](-360/9*4:-.75cm)--node[above,xshift=.6cm,yshift=-.35cm]{\tiny\color{blue}$1^1,2^3$}
(-360/9*5:-.75cm);

\draw[](-360/9*7:-.75cm)--node[below,rotate=60,yshift=.35cm,xshift=-.5cm]{\tiny\color{blue}$1^3,2^2$}
(-360/9*8:-.75cm);

\draw[](-360/9*3:-.75cm)--node[below]{\tiny\color{blue}$3^2$}
(-360/9*3:-.75cm);

\draw[](-360/9*6:-.75cm)--node[above]{\tiny\color{blue}$3^1$}
(-360/9*6:-.75cm);

\draw[](-360/9*9:-.75cm)--node[above,xshift=.3cm,yshift=-.25cm]{\tiny\color{blue}$3^3$}
(-360/9*9:-.75cm);

\foreach \a/\n in {1^1/1,2^1/2,3^1/3,1^3/4,2^3/5,3^3/6,1^2/7,2^2/8,3^2/9}
{
\draw[fill](-360/9*\n:-.75cm)circle(1pt)
node[anchor=-\n*360/9] (\a) {$\a$};
}

\end{tikzpicture},
\qquad
\tikzset{every node/.style={font=\scriptsize}}
\begin{tikzpicture}[baseline=-1mm]

\draw[fill=lightgray](-360/9*1:-.75cm)--
(-360/9*2:-.75cm)--
(-360/9*3:-.75cm)--
(-360/9*1:-.75cm);

\draw[fill=lightgray](-360/9*4:-.75cm)--
(-360/9*5:-.75cm)--
(-360/9*6:-.75cm)--
(-360/9*4:-.75cm);

\draw[fill=lightgray](-360/9*7:-.75cm)--
(-360/9*8:-.75cm)--
(-360/9*9:-.75cm)--
(-360/9*7:-.75cm);

\foreach \a/\n in {1^1/1,2^1/2,3^1/3,1^3/4,2^3/5,3^3/6,1^2/7,2^2/8,3^2/9}
{
\draw[fill](-360/9*\n:-.75cm)circle(1pt)
node[anchor=-\n*360/9] (\a) {$\a$};
}

\end{tikzpicture}
\right)
\]
\begin{tikzpicture}

\draw[fill=lightgray](-360/18*1:-1.5cm)--
(-360/18*2:-1.5cm)--
(-360/18*3:-1.5cm)--
(-360/18*6:-1.5cm)--node[above,rotate=15,xshift=-.15cm,yshift=-.1cm]{\tiny\color{blue}$1^2,2^1$}
(-360/18*1:-1.5cm);

\draw[fill=lightgray](-360/18*7:-1.5cm)--
(-360/18*8:-1.5cm)--
(-360/18*9:-1.5cm)--
(-360/18*12:-1.5cm)--node[below,rotate=80,yshift=.05cm,xshift=.1cm]{\tiny\color{blue} $1^1,2^3$}
(-360/18*7:-1.5cm);

\draw[fill=lightgray](-360/18*13:-1.5cm)--
(-360/18*14:-1.5cm)--
(-360/18*15:-1.5cm)--
(-360/18*18:-1.5cm)--node[below,rotate=-40,yshift=.05cm]{\tiny\color{blue}$1^3,2^2$}
(-360/18*13:-1.5cm);

\draw[fill=lightgray](-360/18*4:-1.5cm)--node[above,yshift=.25cm]{\tiny\color{blue}$3^2$}
(-360/18*5:-1.5cm);

\draw[fill=lightgray](-360/18*10:-1.5cm)--node[below,xshift=.5cm]{\tiny\color{blue}$3^1$}
(-360/18*11:-1.5cm);

\draw[fill=lightgray](-360/18*16:-1.5cm)--node[below,xshift=-.5cm]{\tiny\color{blue}$3^3$}
(-360/18*17:-1.5cm);

\foreach \a/\n in {1^1/1,2^1/2,3^1/3,4^1/4,5^1/5,6^1/6,1^3/7,2^3/8,3^3/9,4^3/10,5^3/11,6^3/12,1^2/13,2^2/14,3^2/15,4^2/16,5^2/17,6^2/18}
{
\draw[fill](-360/18*\n:-1.5cm)circle(1pt)
node[anchor=-\n*360/18] (\a) {$\a$};
}

\end{tikzpicture}
    \caption{An initial $k$-$W$-noncrossing parking function model and its corresponding image under $\nabla$ with $W=G(3,1,3)$ with $k=2$.}
    \label{fig:Image under nabla}
\end{figure}

From~\cite{armstrong2009thesis}, define the \defn{shuffle}, $\shuffle$, of $(\pi_1\leq\dots\leq\pi_k)\in\NC^k(G(d,1,n))$ as follows: Take a model with $kdn$ nodes labelled by creating $k$ copies of each node $i^j$ with a subscript having labels $1,\ldots,k$. We have $k$ nodes for each original node $i^j$ labelled $i^j_1,\ldots,i^j_k$. Connect any nodes with the same subscript $r$ if they are connected in $\pi_r$. To finish, relabel all nodes with the set \[
\{1^1,2^1,\ldots,(kn)^1,1^d,2^d,\ldots,(kn)^d,\ldots,1^2,2^2,\ldots,(kn)^2\}
\]
in clockwise order. We lastly apply the Kreweras complement, $K$, due to~\cite{kreweras1972}. The Kreweras complement applied to a noncrossing partition model begins by taking each node $i^j$ and adding a node labelled $i^{j'}$ preceding $i^j$. Then, connect all possible $i^{j'}$ labelled nodes without crossing any convex hulls existing on the nonprime nodes $i^j$. To finish, delete all nonprime nodes. We have the following lemma due to D.~Armstrong extended for $G(d,1,n)$-noncrossing partitions:

\begin{lemma}[\cite{armstrong2009thesis}]
\label{lemma: k-divisible bij}
    Let $(\pi_1\leq\cdots\leq\pi_k)\in\NC^k(W)$. Then, the partition
    \[
    K\circ \shuffle \circ \Omega^{-1} \circ \delta \circ \Omega (\pi_1\leq\ldots\leq\pi_k)
    \]
    is a $k$-divisible noncrossing partition on $dkn$ nodes and the map $\nabla\coloneq K\circ \shuffle \circ \Omega^{-1} \circ \delta \circ \Omega$ gives a bijection between $k$-multichains in $\NC(W)$ and $k$-divisible noncrossing partitions on 
    \[
    \{1^1,2^1,\ldots,(kn)^1,1^d,2^d,\ldots,(kn)^d,\ldots,1^2,2^2,\ldots,(kn)^2\}.
    \]
    The resulting partition, $\pi$, has the following properties:
    \begin{itemize}
        \item[(1)] Two nodes $i^j,r^l$ belong to the same block in $\pi_1$ if and only if $((i-1)k+1)^j$ and $((r-1)k+1)^l$ belong to the same block in $\pi$. 
        \item[(2)] Given a block $B$ in $\pi_1$ of size $b$ containing $i^j$, the block of $\pi$ containing $((i-1)k+1)^j$ has size $kb$.
    \end{itemize} 
\end{lemma}
The proof of~\Cref{lemma: k-divisible bij} extends easily from that given by D.~Armstrong~\cite{armstrong2009thesis}. We will use the first property in~\Cref{lemma: k-divisible bij} to induce a labelling on $\pi$ from $\pi_1$, thus giving a model for the $k$-$W$-noncrossing parking functions in $G(d,1,n)$. 
\begin{example}
We construct a model from the example in~\Cref{fig:rudimentary fuss model} by applying $\nabla$, whose image is show in~\Cref{fig:Image under nabla}. 
\end{example}

We will refer to the labellings under $\nabla$ by $\ell(B)$ for each block $B\in \pi$. The $(W\times \Z_{kh})$-action corresponds to $W$ permuting the labels, and $\Z_{kh}$ rotating the picture $\frac{2\pi}{kdn}$ clockwise. 

\subsection{Proof of~\Cref{thm: k-bijection} for \texorpdfstring{$G(d,1,n)$}{G(d,1,n)}}
\label{sec:fuss_d1n_bij}
We extend the isomorphism in type $B$ from~\cite{rhoades2014parking} to an isomorphism in the case of $G(d,1,n)$, similarly to our method in~\Cref{sec: d1n bijection}. 

With $G(d,1,n)$ acting on $V=\C^n$ with the standard coordinate functionals ${\bf x}=x_1,\ldots,x_n$, choose 
\[
\Theta_{kh+1}\coloneq(x_1^{kdn+1},\ldots,x_n^{kdn+1}).
\]
By~\Cref{rmk:hsop_inf_fams}, this is an hsop of degree $kdn+1 = kh+1$ carrying $V^*(-kh-1)$. The map $x_i\to \theta_i$ is $W$-equivariant and the subvariety cut out by $(\Theta_{kh+1} - \bf{x})$ is 
\[
V^{\Theta_{kh+1}} \coloneq\{(v_1,\dots,v_n)\in\C^n : v_i = 0 \text{ or } v_i^{dnk} = 1\text{ for all }i\},
\]
with $|V^{\Theta_{kh+1}}|=(dnk+1)^n = (kh+1)^n$. We begin by describing the forward isomorphism  
\[
\Park_\NC^k (W) \xrightarrow{f_k} V^{\Theta_{kh+1}}.
\]

As in the $k=1$ case, we begin by passing our set partition corresponding to the model into an infinite parenthesization and corresponding openers.

\begin{example}
\label{ex:fussd1n_forward_tabulate}
In~\Cref{fig:Image under nabla}, our set partition is 
\[
\{\{1^1,2^1,3^1,6^1\},\{4^1,5^1\},\{1^2,2^2,3^2,6^2\},\{4^2,5^2\},\{1^3,2^3,3^3,6^3\},\{4^3,5^3\}\},
\]
corresponding to the infinite parenthesization
\[
\ldots,(1^1,2^1,3^1,(4^1,5^1),6^1),(1^3,2^3,3^3,(4^3,5^3),6^3),(1^2,2^2,3^2,(4^2,5^2),6^2),\ldots.
\]
From here we tabulate our openers:
\begin{center}
\bgroup
\def\arraystretch{1.5}
\begin{tabular}{|c|c|c|}\hline
nonzero block $B_i$ & opener & $\ell(B_i)$ \\\hline\hline
$\{1^1,2^1,3^1,6^1\}$     & $1^1$ & $\{1^2,2^1\}$     \\\hline
$\{1^2,2^2,3^2,6^2\}$     & $1^2$ & $\{1^3,2^2\}$     \\\hline
$\{1^3,2^3,3^3,6^3\}$     & $1^3$ & $\{1^1,2^3\}$     \\\hline
$\{4^1,5^1\}$ & $4^1$ & $\{3^2\}$    \\\hline
$\{4^2,5^2\}$ & $4^2$ & $\{3^3\}$    \\\hline
$\{4^3,5^3\}$ & $4^3$ & $\{3^1\}$    \\\hline
\end{tabular}.
\egroup
\end{center}
\end{example}

We define our map $f_k$ by passing into the vector $v=(v_1,\ldots,v_n)$ where 
\begin{itemize}
    \item $v_r = 0$ if $r^1\in \ell(B_0)$ where $B_0$ is the zero block, or
    \item $v_r = \lambda^j \omega^i$ if $r^1 \in \ell(B)$ where $B$ is the block with opener $i^j$,
\end{itemize}
with $\lambda = e^{\frac{2\pi \sqrt{-1}}{d}}$ and $\omega = e^{\frac{2\pi \sqrt{-1}}{kn}}$.

\begin{example}
\label{ex:fussd1n_forward_vector}
Continuing with~\Cref{ex:fussd1n_forward_tabulate} we map into the vector
\[
v = (\lambda^3\omega^1, \lambda^1\omega^1, \lambda^3\omega^4).
\]
\end{example}

We now present the inverse map
\[
V^{\Theta_{kh+1}}\xrightarrow{f_k^{-1}}\Park_\NC^k (W).
\]

Given a vector $v \in V^{\Theta_{kh+1}}$, the corresponding infinite parenthesization can be recovered as follows: Let $m_r$ be the number of occurrences of $\lambda^j\omega^r$ for any color $j$ amongst the components of $v$. If $m_r>0$, we have a left parenthesis at every $r^j$ for all $j=1,\ldots,d$. Then, beginning with the smallest nonzero $m_r$, place the corresponding right parenthesis so that the parenthesis pair contains $km_r$ elements, by~\Cref{lemma: k-divisible bij}, excluding any elements already belonging to a parenthesis pair. 

\begin{example}
\label{ex:fussd1n_inverse_partition}
Taking the vector from~\Cref{ex:fussd1n_forward_vector}, we have a left parenthesis to the left of $1^j$ and $4^j$ for an colors $j$. Our multiplicities are 
\begin{equation*}
\begin{split}
    m_1 &= 2 \\
    m_4 &= 1 \\
    m_2 &= m_3=m_5=m_6=0,
\end{split}
\end{equation*}
and with $k=2$, we have the parenthesization:
\[
\ldots,(1^1,2^1,3^1,(4^1,5^1),6^1),(1^3,2^3,3^3,(4^3,5^3),6^3),(1^2,2^2,3^2,(4^2,5^2),6^2),\ldots.
\]
\end{example}

From here we obtain our noncrossing partition and have only the labellings left to recover. We know $r^1 \in w(B)$ if and only if $v_r = i^j$ where $i^j$ is the opener of $B$. In the same way as before, we recover the remaining labellings $r^j$ by the properties of $w\in G(d,1,n)$. 
\begin{example}
\label{ex:fussd1n_inverse_labelling}
Continuing with~\Cref{ex:fussd1n_inverse_partition}, 
\begin{itemize}
    \item $1^1\in \ell(\{1^3,2^3,3^3,6^3\})$,
    \item $2^1\in \ell(\{1^1,2^1,3^1,6^1\})$,
    \item $3^1\in \ell(\{4^3,5^3\})$,
\end{itemize}
and using the properties of $G(d,1,n)$, we obtain the labelling of our model.
\end{example}

Thus the map $\Park^k_\NC (W) \to V^{\Theta_{kh+1}}$ is a bijection. We leave it to the reader to verify this map is $(W\times \Z_{kh})$-equivariant. This concludes the proof of~\Cref{thm: k-bijection} in this case and consequently~\Cref{thm: k-Park is (h+1)^n} by our choice of hsop. 

\section{Type \texorpdfstring{$G(d,d,n)$}{G(d,d,n)}}
\label{sec:gddn}
The infinite family $G(d,d,n)$ is the subfamily of $G(d,1,n)$ with the added condition that the non-zero entries, $h_1,\dots,h_n$ with $h_i\in\mu_d$, in the matrix representation of $w\in G(d,d,n)$ satisfy $h_1\cdots h_n = 1$. In particular, $G(2,2,n)$ recovers the Weyl group of type $D_n$.  As such, we extend isomorphisms in type $D$ from~\cite{armstrong2015parking} to our case. Take Coxeter element $c = [1,2,\ldots,n-1][n]$ with Coxeter number $h = (n-1)d$. 

\subsection{Visualizing type \texorpdfstring{$G(d,d,n)$}{G(d,d,n)} noncrossing partitions}
\label{sec:gddn_model}
We describe the model for $G(d,d,n)$-noncrossing partitions due to D.~Bessis and R.~Corran~\cite{bessis2006eer}. Let $X$ be a finite subset of $\C$. 
\begin{definition}
    A partition $\pi$ of $X$ is \defn{noncrossing} if for any $B,B'\in \pi$, if the intersection of the convex hulls of $B$ and $B'$ are nonempty, then $B = B'$.
\end{definition}
Let $\operatorname{NCP}_X$ be the set of all noncrossing partitions of a finite subset $X\subset \C$. 
\begin{definition}[\cite{bessis2006eer}, Definition 1.15]
    Let $\pi^{\flat}$ be obtained from $\pi$ by removing 0. Then, 
    \[
    \NC(G(d,d,n))\coloneq\{\pi\in\operatorname{NCP}_{\mu_{d(n-1)}\cup\{0\}}:\pi^{\flat}\in \NC(G(d,1,n-1))\}.
    \]
\end{definition}

D.~Bessis and R.~Corran note in~\cite{bessis2006eer} that their definition is equivalent to that due to C.~Athanasiadis V.~Reiner for type $D_n$ in~\cite{reiner2004Dn}. There is a map illustrating this equivalence, which in general gives us a labelled model for each $G(d,d,n)$-noncrossing partition. We will illustrate this mapping and use the labelled model to enumerate in $G(d,d,n)$.

For each $d$ and $n$, fix a $z_0 \in \mu_{d(n-1)}$. Beginning at this $z_0$, label the nodes belonging to $\mu_{d(n-1)}$ in clockwise order by 
\[
1^1,\ldots,(n-1)^1,1^d,\ldots,(n-1)^d,1^{d-1},\ldots,(n-1)^{d-1},\ldots,1^2,\ldots,(n-1)^2.
\]
Then add $n^1$ to the unique block $B^*$ containing ${0}$. Finally, add $n^j$ to the block $B=\operatorname{rot}_j(B^*)$ where $\operatorname{rot}_j$ is counterclockwise rotation by $\frac{2\pi j}{(n-1)}$ for $j=1,2,\ldots,d-1$. This gives a labelled model for noncrossing partitions in $G(d,d,n)$ that is symmetric under $\frac{2\pi}{(n-1)}$ rotation. For the inverse, remove the central node from any block $B$ containing $n^j$ with $j\neq 1$ and not containing $n^1$. After, remove the labels and we recover the unlabelled noncrossing partition. 

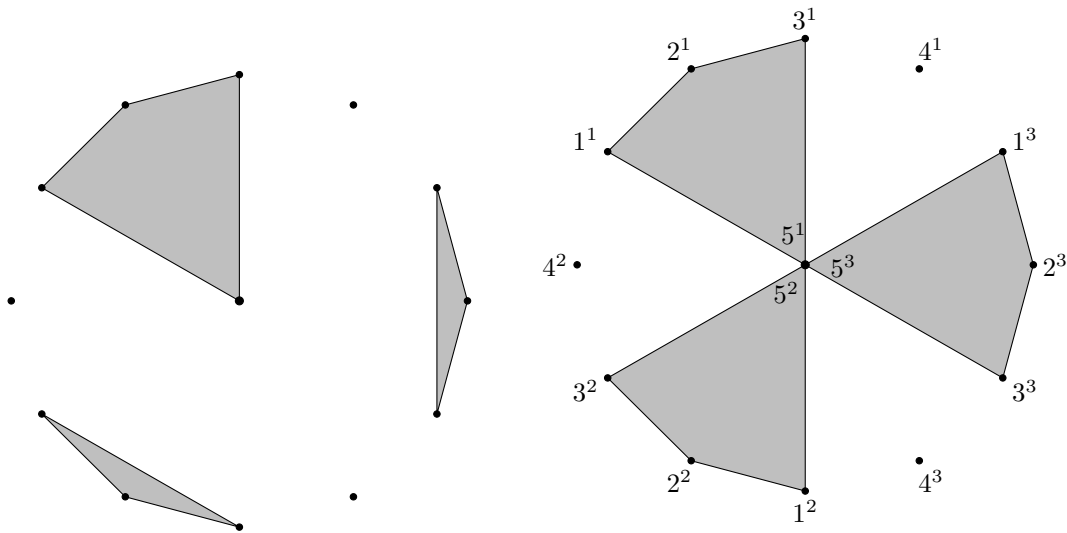
\begin{figure}[h]
    \centering
        \begin{tikzpicture}
\draw[fill=lightgray](-360/12*1:-3cm)--
(-360/12*2:-3cm)--
(-360/12*3:-3cm)--
(0,0)--
(-360/12*1:-3cm);

\draw[fill=lightgray] (-360/12*5:-3cm)--
(-360/12*6:-3cm) --
(-360/12*7:-3cm)--
(-360/12*5:-3cm);

\draw[fill=lightgray] (-360/12*9:-3cm)--
(-360/12*10:-3cm) --
(-360/12*11:-3cm)--
(-360/12*9:-3cm);

\node[circle,draw=black, fill=black, inner sep=0pt,minimum size=3pt] (0) at (0,0) {};

\foreach \a/\n in {1^1/1,2^1/2,3^1/3,4^1/4,1^3/5,2^3/6,3^3/7,4^3/8,1^2/9,2^2/10,3^2/11,4^2/12}
{
\draw[fill](-360/12*\n:-3cm)circle(1.2pt);
}

\end{tikzpicture}
\qquad
        \begin{tikzpicture}
\draw[fill=lightgray] (-360/12*1:-3cm)--
(-360/12*2:-3cm) --
(-360/12*3:-3cm)--
(0,0)node[above,xshift=-.15cm,yshift=.15cm]{$5^1$}--
(-360/12*1:-3cm);

\draw[fill=lightgray] (-360/12*5:-3cm)--
(-360/12*6:-3cm) --
(-360/12*7:-3cm)--
(0,0)node[below,xshift=.5cm,yshift=.25cm]{$5^3$}--
(-360/12*5:-3cm);

\draw[fill=lightgray] (-360/12*9:-3cm)--
(-360/12*10:-3cm) --
(-360/12*11:-3cm)--
(0,0)node[below,xshift=-.25cm,yshift=-.1cm]{$5^2$}--
(-360/12*9:-3cm);

\node[circle,draw=black, fill=black, inner sep=0pt,minimum size=3pt] (0) at (0,0) {};

\foreach \a/\n in {1^1/1,2^1/2,3^1/3,4^1/4,1^3/5,2^3/6,3^3/7,4^3/8,1^2/9,2^2/10,3^2/11,4^2/12}
{
\draw[fill](-360/12*\n:-3cm)circle(1.2pt)
node[anchor=-\n*360/12] (\a) {$\a$};
}

\end{tikzpicture}
    \caption{Unlabelled (left) to labelled (right) noncrossing partition example for $G(3,3,5)$.}
    \label{fig:ncp example G335}
\end{figure}

\begin{example}
    Consider~\Cref{fig:ncp example G335}, which shows an unlabelled noncrossing partition and its corresponding labelled noncrossing partition.
\end{example}

From here, we label the noncrossing partition model with each block's image under a given $w$ to pictorially represent a $G(d,d,n)$-noncrossing parking function. The $(W\times C)$-action in this case corresponds to $W$ permuting the labels of each block of $\pi$ and $C$ rotating the outer boundary clockwise and the inner $\{n^j\}$ labels counterclockwise.

\begin{example}
Consider for example,
\begin{equation}
\label{eqn:w_for_g335}
   w= \left(
\begin{matrix}
1^1 & 2^1 & 3^1 & 4^1 & 5^1 \\
2^2 & 3^1 & 5^1 & 1^3 & 4^1
\end{matrix}
\right).
\end{equation}
The $G(d,d,n)$-noncrossing parking functions fall into one of three cases:
\begin{itemize}
    \item $\pi$ does not contain a zero block nor singletons $\{n^j\}$,
    \item $\pi$ contains a zero block, or
    \item $\pi$ contains singletons $\{n^j\}$.
\end{itemize}
Examples of these three cases are shown in~\Cref{fig:ncpf_example_g335}. In these models we may exclude the labels $\{n^j\}$ if a zero block is present as $n^j$ must belong to the zero block if it exists. Similarly, if the singletons $\{n^j\}\in \pi$, we only label the central node by $n^1$.
\end{example}
    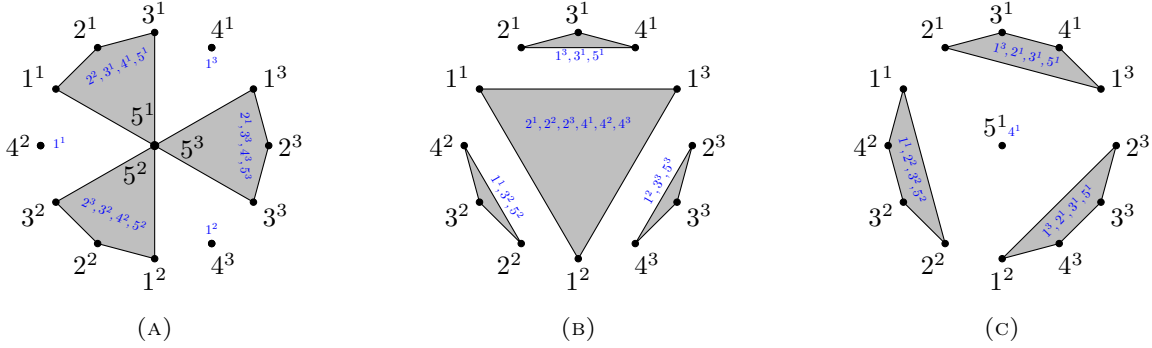
\begin{figure}[h]
        \centering
        \begin{subfigure}[h]{0.3\textwidth}
            \centering
            \begin{tikzpicture}

\draw[fill=lightgray] (-360/12*1:-1.5cm)--
node[below,rotate=25,xshift=.5cm]{\scalebox{.5}{\color{blue}$2^2,3^1,4^1,5^1$}}
(-360/12*2:-1.5cm) --
(-360/12*3:-1.5cm)--
(0,0)node[above,xshift=-.15cm,yshift=.15cm]{$5^1$}--
(-360/12*1:-1.5cm);

\draw[fill=lightgray] (-360/12*5:-1.5cm)--
node[below,rotate=-90,xshift=.35cm]{\scalebox{.5}{ \color{blue}$2^1,3^3,4^3,5^3$}}
(-360/12*6:-1.5cm) --
(-360/12*7:-1.5cm)--
(0,0)node[below,xshift=.5cm,yshift=.25cm]{$5^3$}--
(-360/12*5:-1.5cm);

\draw[fill=lightgray] (-360/12*9:-1.5cm)--
node[above,rotate=-25,xshift=-.35cm,yshift=.15cm]{\scalebox{.5}{\color{blue}$2^3,3^2,4^2,5^2$}}
(-360/12*10:-1.5cm) --
(-360/12*11:-1.5cm)--
(0,0)node[below,xshift=-.25cm,yshift=-.1cm]{$5^2$}--
(-360/12*9:-1.5cm);

\draw[fill](-360/12*4:-1.5cm)circle(1.2pt)node[below]{\scalebox{.5}{\color{blue}$1^3$}};
\draw[fill](-360/12*8:-1.5cm)circle(1.2pt)node[above]{\scalebox{.5}{\color{blue}$1^2$}};
\draw[fill](-360/12*12:-1.5cm)circle(1.2pt)node[above,xshift=.25cm,yshift=-.15cm]{\scalebox{.5}{\color{blue}$1^1$}};

\node[circle,draw=black, fill=black, inner sep=0pt,minimum size=3pt] (0) at (0,0) {};

\foreach \a/\n in {1^1/1,2^1/2,3^1/3,4^1/4,1^3/5,2^3/6,3^3/7,4^3/8,1^2/9,2^2/10,3^2/11,4^2/12}
{
\draw[fill](-360/12*\n:-1.5cm)circle(1.2pt)
node[anchor=-\n*360/12] (\a) {$\a$};
}

\end{tikzpicture}
            \caption{}
            \label{subfig: ncpf ex subfig A}
        \end{subfigure}\hfill
        \begin{subfigure}[h]{0.3\textwidth}
            \centering
            \begin{tikzpicture}
\node[circle,draw=black, fill=black, inner sep=0pt,minimum size=1pt] (0) at (0,0) {};

\draw[fill=lightgray] (-360/12*1:-1.5cm)--
node[below,yshift=-.25cm]{\scalebox{.5}{\color{blue}$2^1,2^2,2^3,4^1,4^2,4^3$}}
(-360/12*5:-1.5cm)--
(-360/12*9:-1.5cm)--
(-360/12*1:-1.5cm);

\draw[fill=lightgray] (-360/12*2:-1.5cm)--
node[below,xshift=.4cm]{\scalebox{.5}{\color{blue}$1^3,3^1,5^1$}}
(-360/12*3:-1.5cm) --
(-360/12*4:-1.5cm)--
(-360/12*2:-1.5cm);

\draw[fill=lightgray] (-360/12*6:-1.5cm)--
node[below,rotate=60,xshift=-.25cm,yshift=.45cm]{\scalebox{.5}{\color{blue}$1^2,3^3,5^3$}}
(-360/12*7:-1.5cm) --
(-360/12*8:-1.5cm)--
(-360/12*6:-1.5cm);

\draw[fill=lightgray] (-360/12*10:-1.5cm)--
node[below,rotate=-60,xshift=-.25cm,yshift=.45cm]{\scalebox{.5}{\color{blue}$1^1,3^2,5^2$}}
(-360/12*11:-1.5cm) --
(-360/12*12:-1.5cm)--
(-360/12*10:-1.5cm);

\foreach \a/\n in {1^1/1,2^1/2,3^1/3,4^1/4,1^3/5,2^3/6,3^3/7,4^3/8,1^2/9,2^2/10,3^2/11,4^2/12}
{
\draw[fill](-360/12*\n:-1.5cm)circle(1.2pt)
node[anchor=-\n*360/12] (\a) {$\a$};
}
\end{tikzpicture}
            \caption{}
            \label{subfig: ncpf ex subfig B}
        \end{subfigure}\hfill
        \begin{subfigure}[h]{0.3\textwidth}
            \centering
            \begin{tikzpicture}
\draw[fill](0,0)circle(1.2pt)
node[above] (0) {$5^1$\scalebox{.5}{\color{blue}$4^1$}};

\draw[fill=lightgray] (-360/12*2:-1.5cm)--
(-360/12*3:-1.5cm)--node[below,rotate=-15]{\scalebox{.5}{\color{blue}$1^3,2^1,3^1,5^1$}}
(-360/12*4:-1.5cm)--
(-360/12*5:-1.5cm)--
(-360/12*2:-1.5cm);

\draw[fill=lightgray] (-360/12*6:-1.5cm)--
(-360/12*7:-1.5cm)--node[above,rotate=45]{\scalebox{.5}{\color{blue}$1^3,2^1,3^1,5^1$}}
(-360/12*8:-1.5cm)--
(-360/12*9:-1.5cm)--
(-360/12*6:-1.5cm);

\draw[fill=lightgray] (-360/12*10:-1.5cm)--
(-360/12*11:-1.5cm)--node[above,rotate=-75]{\scalebox{.5}{\color{blue}$1^1,2^2,3^2,5^2$}}
(-360/12*12:-1.5cm)--
(-360/12*1:-1.5cm)--
(-360/12*10:-1.5cm);

\foreach \a/\n in {1^1/1,2^1/2,3^1/3,4^1/4,1^3/5,2^3/6,3^3/7,4^3/8,1^2/9,2^2/10,3^2/11,4^2/12}
{
\draw[fill](-360/12*\n:-1.5cm)circle(1.2pt)
node[anchor=-\n*360/12] (\a) {$\a$};
}

\end{tikzpicture}
            \caption{}
            \label{subfig: ncpf ex subfig C}
        \end{subfigure}
        \caption{Three noncrossing parking functions in $G(3,3,5)$ with $w$ as in \Cref{eqn:w_for_g335}}.
        \label{fig:ncpf_example_g335}
    \end{figure}

\subsection{Proof of \Cref{thm: bijection} for \texorpdfstring{$G(d,d,n)$}{G(d,d,n)}}
\label{sec:ddn bijection}
We follow the treatment of type $D_n$ by~\cite{armstrong2015parking} and extend it to our models in $G(d,d,n)$. Let $W=G(d,d,n)$ act on $V=\C^n$ with the standard coordinate functionals ${\bf x}=x_1,\ldots,x_n$ and choose
\[
(\Theta_{h+1}) = (\theta_1,\dots,\theta_n)\coloneq(x_1^{d(n-1)+1},\dots,x_n^{d(n-1)+1}).
\]
Again by~\Cref{rmk:hsop_inf_fams}, this is an hsop of degree $d(n-1)+1 = h+1$ carrying $V^*(-h-1)$ and the subvariety cut out by $(\Theta_{h+1}-\bf{x})$,
\[
V^{\Theta_{h+1}}\coloneq\{(v_1,\ldots,v_n)\in\C^n:v_i=0\text{ or }v_i^{d(n-1)}=1\},
\]
contains $(d(n-1)+1)^n=(h+1)^n$ distinct elements. We describe first the forward isomorphism 
\[
\Park_\NC^k (W) \xrightarrow{f_d} V^{\Theta_{h+1}}.
\]

Given a noncrossing parking function $[w,\pi]$, first take the underlying noncrossing partition and exclude the elements $n^j$ for any color $j=1,\ldots,d$. Then, map the resulting partition into an infinite parenthesization of the string:
\[
\dots 1^1,\dots,(n-1)^1,1^d,\dots,(n-1)^d, 1^{d-1},\dots,(n-1)^{d-1},\dots,1^2,\dots,(n-1)^2,1^1\dots.
\]
Using the same definition of opener as before, we tabulate the nonzero blocks, their openers, and their images under $w$.

\begin{example}
\label{ex:ddn_forward_tabulate}
Consider each of the three examples in~\Cref{fig:ncpf_example_g335}. The exponent of the block $B_i$ indicates the example for which the column is taken from. We tabulate the nonzero block openers as follows:
\begin{center}
\bgroup
\def\arraystretch{1.5}
\begin{tabular}{|c|c|c||c|c|c|}\hline
$B_i^\text{A}$ & opener & $w(B_i^\text{A})$ & $B_i^\text{B}$ & opener & $w(B_i^\text{B})$ \\\hline\hline
\small $\{1^1,2^1,3^1\}$ & \small $1^1$ & \small$\{2^2,3^1,4^1,5^1\}$  &\small $\{2^1,3^1,4^1\}$ &\small $2^1$ &\small $\{1^3,3^1,5^1\}$ \\\hline
\small $\{1^2,2^2,3^2\}$ & \small $1^2$ & \small$\{2^3,3^2,4^2,5^2\}$   &\small $\{2^2,3^2,4^2\}$ &\small $2^2$ &\small $\{1^1,3^2,5^2\}$ \\\hline
\small $\{1^3,2^3,3^3\}$ & \small $1^3$ & \small$\{2^1,3^3,4^3,5^3\}$   &\small $\{2^3,3^3,4^3\}$ &\small $2^3$ &\small $\{1^2,3^3,5^3\}$  \\\hline
\small $\{4^1\}$ & \small $4^1$ & \small$\{1^3\}$   & & & \\\hline
\small $\{4^2\}$ & \small $4^2$ & \small$\{1^1\}$   & & & \\\hline
\small $\{4^3\}$ & \small $4^3$ & \small$\{1^2\}$   & & & \\\hline
\end{tabular}
\egroup
\end{center}
\begin{center}
\bgroup
\def\arraystretch{1.5}
\begin{tabular}{|c|c|c|}\hline
$B_i^\text{C}$ & opener & $w(B_i^\text{C})$\\\hline\hline
\small $\{1^3,2^1,3^1,4^1\}$&\small $2^1$ &\small $\{1^3,2^1,3^1,5^1\}$ \\\hline
\small $\{1^1,2^2,3^2,4^2\}$&\small $2^2$ &\small $\{1^1,2^2,3^2,5^2\}$  \\\hline
\small $\{1^2,2^3,3^3,4^3\}$&\small $2^3$ &\small $\{1^2,2^3,3^3,5^3\}$  \\\hline
\small $\{5^1\}$&\small $5^1$ &\small $\{4^1\}$  \\\hline
\small $\{5^2\}$&\small $5^2$ &\small $\{4^2\}$  \\\hline
\small $\{5^3\}$&\small $5^3$ &\small $\{4^3\}$  \\\hline
\end{tabular}
\egroup
\end{center}
\end{example}

We define the map $f_d$ by passing into the vector $v=(v_1,\ldots,v_n)$ where
\begin{itemize}
    \item $v_r=0$ if $r^1\in w(B_0)$ where $B_0$ is the zero block,
    \item $v_r = 0$ if $r^1 \in w(\{n^j\})$ for any color $j=1,\ldots,d$, or
    \item $v_r = \lambda^j\omega^i$ if $r^1\in w(B)$ where $B$ is a nonzero block opened by $i^j$,
\end{itemize}
with $\lambda\coloneq e^{\frac{2\pi \sqrt{-1}}{d}}$ and $\omega \coloneq e^{\frac{2\pi \sqrt{-1}}{n-1}}$.
\begin{example}
\label{ex:ddn_forward_vector}
Continuing with~\Cref{ex:ddn_forward_tabulate}, we have the following images of our three examples from~\Cref{fig:ncpf_example_g335}:
\begin{itemize}
    \item[(\ref{subfig: ncpf ex subfig A})] $v = (\lambda^2\omega^4,\lambda^3\omega^1,\lambda^1\omega^1,\lambda^1\omega^1,\lambda^1\omega^1)$,
    \item[(\ref{subfig: ncpf ex subfig B})]$v = (\lambda^2\omega^2,0,\lambda^1\omega^2,0,\lambda^1\omega^2)$, \hspace{.15cm}and 
    \item[(\ref{subfig: ncpf ex subfig C})]$v = (\lambda^2\omega^2,\lambda^1\omega^2,\lambda^1\omega^2,0,\lambda^1\omega^2)$.
\end{itemize}
\end{example}

To show this map is indeed a bijection, we construct its inverse
\[
V^{\Theta_{h+1}} \xrightarrow{f_d^{-1}} \Park_\NC^k (W).
\]

Given a vector $v\in V^{\Theta_{h+1}}$, let $m_r$ be the number of occurrences of $\lambda^j\omega^r$ among the components of $v$ for any color $j=1,\ldots,d$. Let $z$ be the number of zero components of $v$. We have three cases to consider for the inverse:
\begin{itemize}
    \item[\bf{Case 1:}] $z=0$.
    In this case, the noncrossing parking function must not have a nonzero block and $n^j$ must not be a singleton for any color $j=1,\ldots,d$. Then, our multiplicities sum to $n$ and for each $r$ with $m_r\geq 2$, the adjusted multiplicities $(m_1,m_2,\ldots,m_r-1,\ldots,m_{n-1})$ give a noncrossing partition for type $G(d,1,n-1)$ with its labellings. One can check that there is exactly one choice, $r_*$, among the multiplicities $m_r\geq 2$ such that when the central vertex is added to the block opened by $r_*^j$ for each color $j$, the picture remains noncrossing for $G(d,d,n)$. The choice of which $n^\ell$ to add to the block opened by $r_*^1$ is determined by the fact that the entries of $w$ in its matrix representation must multiply to 1. 
    \begin{example}
    \label{ex:ddn_inverse_exampleA}
    The example in~\Cref{subfig: ncpf ex subfig A} falls under this case. We have multiplicities 
    \begin{equation*}
`   \begin{split}
    m_1 &= 4 \\
    m_4 &= 1 \\
    m_2 &= m_3=0.
    \end{split}
    \end{equation*}
    Which give us the noncrossing parking function in $G(d,1,n-1)$ from multiplicities $(m_1-1,m_2,m_3,m_4)=(3,0,0,1)$ and labellings where the blocks opened by $1^j$ have one more label than the block size. We add $5^1$ to the block $\{1^1,2^1,3^1\}$ as any other choice would violate the condition on $w$.
    \end{example}
    
    \item[\bf{Case 2:}] $z\geq2$. In this case, the corresponding noncrossing parking function must have a zero block of size $dz$. We use our multiplicities in the same way as in $G(d,1,n-1)$ to recover the noncrossing partition and its labelling. To recover the corresponding $G(d,d,n)$ noncrossing parking function we simply add each $n^j$ to the zero block for $j=1,\ldots,d$. The zero block will now be labelled by $d$ fewer elements and we may add the missing element in all of its colors to the label of the zero block to recover our noncrossing parking function. 
    \begin{example}
        The example in~\Cref{subfig: ncpf ex subfig B} falls into this case. We reconstruct the underlying $G(3,1,4)$-noncrossing partition in the same way along with its labellings. We simply add the singletons $\{5^j\}$ to the zero block and the labellings that are unused to the label of the zero block.
    \end{example}
    
    \item[\bf{Case 3:}] $z=1$. In this case, we have $d$ singletons $\{n^j\}$ for each $j=1,\dots,d$. Our multiplicities can then be treated exactly as in the case of $G(d,1,n-1)$ recovering a noncrossing parking function in this type. We then determine the labelling of the singleton $\{n^1\}$ by the fact that the entries in the matrix representation of $w$ must multiply to equal 1. This induces the labelling for each $\{n^j\}$ for $j=2,\ldots, d$.
    \begin{example}
    The example in~\Cref{subfig: ncpf ex subfig C} falls into this case. We first reconstruct the underlying $G(3,1,4)$-noncrossing parking function and finish by labelling the singletons. We see that $w(5^1) = 4^1$ as otherwise the product of the entries in the matrix representation of $w$ would not equal 1. 
    \end{example}
\end{itemize}

Thereby, we have shown the inverse of this map and hence that it is a bijection. We leave it to the reader to verify this map is $(W\times C)$-equivariant. Thus, we have proven~\Cref{thm: bijection} in this case and~\Cref{thm: Park is (h+1)^n} by our choice of hsop.

\section{Fuss Analogues of Type \texorpdfstring{$G(d,d,n)$}{G(d,d,n)}}
\label{sec:fuss_gddn}
As with our treatment of the Fuss case for $G(d,1,n)$ we extend the method's used by B.~Rhoades in~\cite{rhoades2014parking} in type $D$ to $G(d,d,n)$. 

\subsection{Visualizing Fuss \texorpdfstring{$G(d,d,n)$}{G(d,d,n)} Noncrossing Partitions}
We begin with a naive model giving $k$ noncrossing parking functions corresponding to a multichain $(\pi_1\leq\cdots\leq\pi_k)$, with the first labelled by its image under $w\in G(d,d,n)$. To begin, pass $(\pi_1\leq\cdots,\leq\pi_k)$ to $(w_1\leq\cdots\leq w_k)$ under $\Omega$ as defined in~\Cref{ex:gd1n_wxc_action}. We extend a map due to~\cite{KM2010}\cite{Kim2011} that B.~Rhoades uses in handling type $D$ (see~\cite{rhoades2014parking}, Theorem 7.1).

Let $k\geq 1$. Given any $u\in G(d,d,n)$ and any $1\leq m \leq k$, define $\tau_{k,m}(u)\in G(d,d,kn)$ by letting the cycles of $\tau_{k,m}(u)$ be obtained by from the cycles of $u$ by replacing each $i^j$ by $((i-1)k+m)^j$. For any $(w_1\leq\cdots\leq w_k)\in\NC^k(W)$ define
\begin{equation}\label{tau map}
\begin{split}
    \tau(w_1\leq\cdots\leq w_k)&\coloneq
    [1,2,\ldots,k(n-1)][k(n-1)+1,k(n-1)+2,\ldots,kn]\\
    &\cdot \tau_{k,1}(w_1^{-1}w_2)^{-1}\cdots \tau_{k,k-1}(w_{k-1}^{-1}w_k)^{-1}\tau_{k,k}(w_k^{-1}c)^{-1}.
\end{split}
\end{equation}

We may then visualize the cycles of $\tau(w_1\leq\cdots\leq w_k)$ as a noncrossing partition on an annulus as follows: First, take an annulus with $dk(n-1)$ outer vertices and $dk$ inner vertices. Label the outer vertices clockwise by
\[
1^1,2^1,\ldots,(k(n-1))^1,1^d,2^d,\ldots,(k(n-1))^d,1^{d-1},2^{d-1},\ldots,(k(n-1))^{d-1},\ldots,1^2,2^2,\ldots,(k(n-1))^2,
\]
and label the inner annulus counterclockwise with 
\[
(k(n-1)+1)^1,\ldots,(kn)^1,(k(n-1)+1)^2,\ldots,(kn)^2,\ldots,(k(n-1)+1)^d,\ldots,(kn)^d.
\]

Let $\pi$ be a partition of $[kn]^d$. 

\begin{definition}[\cite{rhoades2014parking}\cite{KM2010}\cite{Kim2011}]\label{def: kn-noncrossing partition}
We say a block of $B$ of $\pi$ is \defn{inner (outer)} if it consists only of vertices on the inner (outer) boundary. We say a block $B$ of $\pi$ is annular if it is not inner, outer, or the zero block. The partition $\pi$ said to be \defn{$kn$-noncrossing} if the following conditions hold:
\begin{enumerate}
    \item If $B\in \pi$, then $B^r\coloneq\{i^{j+r\mod d}:i^j\in B\}\in\pi$ for all $r=1,\ldots,d$.
    \item The convex hulls of the vertices on the annulus corresponding to distinct blocks of $\pi$ do not intersect.
    \item If $\pi$ contains a zero block, $B_0$, then $B_0$ must contain all vertices on the inner boundary along with $dk$ additional vertices on the outer boundary.
    \item When the vertices of any block of $\pi$ are read in cyclic order from the embedding into the annulus, the consecutive vertices number (ignoring their color) represent consecutive residue classes modulo $k$.
    \item If every block of $\pi$ is either inner or outer, then the outer blocks of $\pi$ determines the inner blocks in the following way: We call an outer block \defn{visible} if it can be connected to the inner boundary without crossing any other outer blocks. Let $B$ be a visible block and let $b^r\in B$ be the last vertex read in clockwise order. Then the $d$ inner blocks of $\pi$ are the unique blocks containing consecutive vertices on the inner annulus each of size $k$ ending with $a^j$ for each $j=1,\dots,d$ where $a\equiv b \mod k$.
\end{enumerate}
\end{definition}

With this definition, we have a bijection between the $k$-$W$-noncrossing partitions and these $kn$-noncrossing partitions.  

\begin{theorem}[\cite{rhoades2014parking}\cite{KM2010}\cite{Kim2011}]\label{Thm: tau bijection}
    With $W=G(d,d,n)$, there is a bijection from the set $\NC^k(W)$ of $k$-$W$-noncrossing partitions to the set of $kn$-noncrossing partitions as defined above by sending the multichain $(\pi_1\leq\cdots\leq\pi_k)$ to its corresponding multichain in cycles $(w_1\leq\cdots\leq w_k)$ and allowing $\pi$ to be the partition whose blocks correspond to the cycles of $\tau(w_1\leq\cdots\leq w_k)$. Further, we obtain information about $\pi$ given $\pi_1$:
    \begin{itemize}
        \item $\pi$ has a zero block if and only if $\pi_1$ does. If the zero block is present in $\pi$, it has size $k$ times the size of the zero block of $\pi_1$.
        \item $\pi$ has inner blocks if and only if $\{n^j\}$ are singletons for each $j=1,\dots,d$ in $\pi_1$.
        \item If $\pi_1$ has $r_i$ blocks of size $i$ then $\pi$ has $r_i$ blocks of size $ki$.
    \end{itemize}
\end{theorem}

We now extend this definition of a $k$-$W$-noncrossing partition to give us a model for our \defn{$k$-$W$-noncrossing parking function} by extending the work of B.~Rhoades~\cite[Lemma 7.3]{rhoades2014parking}.

\begin{lemma}
[\cite{rhoades2014parking}, Lemma 7.3]\label{lemma:mknd_conditions}
    Let $W=G(d,d,n)$ and $k\geq1$. Define the set $M(k,n,d)$ of \defn{$kn$-noncrossing parking functions} as the set of pairs $(\pi,\ell)$ where $\pi$ is a $kd$-noncrossing partition and $\ell$ is a chosen labelling of the blocks in $\pi$ with the subset $\{1^1,\ldots,n^1,1^2,\ldots,n^2,\ldots,1^d,\ldots,n^d\}$ satisfying
    \begin{enumerate}
        \item For any $B\in \pi$, $|B|=k|\ell(B)|$.
        \item $\{1^1,\ldots,n^1,\ldots,1^d,\ldots,n^d\} = \cup_{B\in\pi}\ell(B)$.
        \item For any $B\in \pi$, $\ell(B^r)=\ell(B)^r$ where $B^r\coloneq\{i^{j+r}:i^j\in B\}$.
        \item Let $Y\coloneq\{i^j:i=1,\ldots,kn, j=1,\ldots,d, \text{ and } i\equiv 1 \mod k\}$. There exists a bijection $\phi:Y\to\{1^1,\ldots,n^1,\ldots,1^d,\ldots,n^d\}$ such that $\phi(i^j) = \phi(i^1)^j$ and $\phi$ restricts to a bijection on $B\cap Y \to \ell(B)$ for all blocks $B\in \pi$. Call \defn{an increase in color by $r$} when $\phi(i^j) = i^{j+r \mod d}$ for some $i^j$. Then the total increases in color of $\phi$ being the sum of all such $r$'s is divisible by $d$. 
    \end{enumerate}
    The $(W\times \Z_{kh})$-action on $M(k,n,d)$ is given by letting $W$ act on the labels and $\Z_{kh}$ act by clockwise rotation on the outer boundary, and counterclockwise rotation on the inner boundary. There exists a $(W\times \Z_{kh})$-equivariant bijection $M(k,n,d)\cong \Park^k_\NC (W)$.
\end{lemma} 

We have the following cases for all partitions $\pi$ given the above definition. If $\pi$ is $kn$-noncrossing, then $\pi$ falls into one of the following three cases:
\begin{itemize}
    \item $\pi$ contains a zero block and all other blocks of $\pi$ are inner,
    \item $\pi$ contains inner blocks and all other blocks of $\pi$ are outer, or
    \item $\pi$ contains annular blocks and all other blocks of $\pi$ are outer.
\end{itemize}

The proofs of~\Cref{Thm: tau bijection} and~\Cref{lemma:mknd_conditions}, due to~\cite{KM2010,Kim2011} and~\cite{rhoades2014parking} respectively, extend without difficulty so we omit them here. 

 \begin{example}
     In~\Cref{fig:3_mknd_examples}, we see three elements belonging to $M(3,4,3)$ each falling into a separate case. From left to right,~\Cref{subfig:3mknd_a} has only annular blocks,~\Cref{subfig:3mknd_b} has only outer or inner blocks, and~\Cref{subfig:3mknd_c} has a zero block.
 \end{example}
    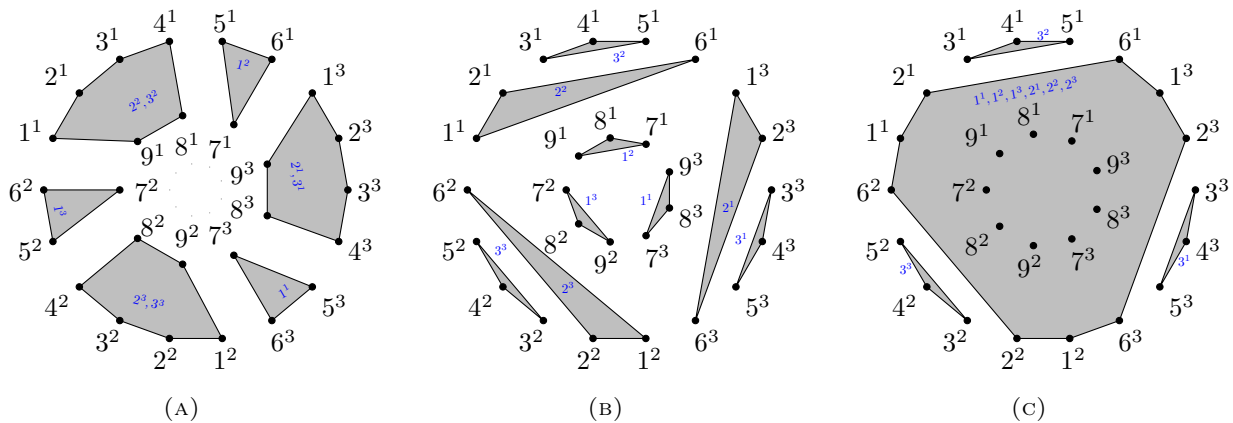
\begin{figure}[h]
        \centering
        \begin{subfigure}[h]{0.3\textwidth}
            \centering
            \begin{tikzpicture}
             \draw[fill=lightgray] (-360/18*1:-2cm)--
    (-360/18*2:-2cm)--node[below,xshift=.5cm,yshift=-.15cm,rotate=25]{\scalebox{.5}{\color{blue}$2^2,3^2$}}
    (-360/18*3:-2cm)--
    (-360/18*4:-2cm)--
    (-360/9*2:-1cm)--
    (-360/9*1:-1cm)--
    (-360/18*1:-2cm);

    \draw[fill=lightgray] (-360/18*7:-2cm)--
    (-360/18*8:-2cm)--node[below,xshift=-.4cm,yshift=-.15cm,rotate=-85]{\scalebox{.5}{\color{blue}$2^1,3^1$}}
    (-360/18*9:-2cm)--
    (-360/18*10:-2cm)--
    (-360/9*5:-1cm)--
    (-360/9*4:-1cm)--
    (-360/18*7:-2cm);

    \draw[fill=lightgray] (-360/18*13:-2cm)--
    (-360/18*14:-2cm)--node[above,yshift=.15cm,rotate=-15]{\scalebox{.5}{\color{blue}$2^3,3^3$}}
    (-360/18*15:-2cm)--
    (-360/18*16:-2cm)--
    (-360/9*8:-1cm)--
    (-360/9*7:-1cm)--
    (-360/18*13:-2cm);

    \draw[fill=lightgray] (-360/18*5:-2cm)--
    node[below,rotate=-15]{\scalebox{.5}{\color{blue}$1^2$}}
    (-360/18*6:-2cm)--
    (-360/9*3:-1cm)--
    (-360/18*5:-2cm);

    \draw[fill=lightgray] (-360/18*11:-2cm)--
    node[above,rotate=30]{\scalebox{.5}{\color{blue}$1^1$}}
    (-360/18*12:-2cm)--
    (-360/9*6:-1cm)--
    (-360/18*11:-2cm);

    \draw[fill=lightgray] (-360/18*17:-2cm)--node[below,xshift=.35cm,yshift=.1cm,rotate=-80]{\scalebox{.5}{\color{blue}$1^3$}}
    (-360/18*18:-2cm)--
    (-360/9*9:-1cm)--
    (-360/18*17:-2cm);
    
   \foreach \a/\n in {1^1/1,2^1/2,3^1/3,4^1/4,5^1/5,6^1/6,1^3/7,2^3/8,3^3/9,4^3/10,5^3/11,6^3/12,1^2/13,2^2/14,3^2/15,4^2/16,5^2/17,6^2/18}
{
\draw[fill](-360/18*\n:-2cm)circle(1.2pt)
node[anchor=-\n*360/18] (\a) {$\a$};
}

   \foreach \a/\n in {9^1/1,8^1/2,7^1/3,9^3/4,8^3/5,7^3/6,9^2/7,8^2/8,7^2/9}
{
\draw[fill](-360/9*\n:-1cm)circle(1.2pt)
node[anchor=-\n*360/9] (\a) {};
}
 \foreach \a/\n in {9^1/1,8^1/2,7^1/3,9^3/4,8^3/5,7^3/6,9^2/7,8^2/8,7^2/9}
{
\draw[fill](-360/9*\n:-.35cm)circle(0pt)
node[anchor=-\n*360/9] (\a+100) {$\a$};
}
\end{tikzpicture}
            \caption{}
            \label{subfig:3mknd_a}
        \end{subfigure}\hfill
        \begin{subfigure}[h]{0.3\textwidth}
            \centering
            \begin{tikzpicture}
        \draw[fill=lightgray] (-360/9*1:-.7cm)--
    (-360/9*2:-.7cm)--
    node[below]{\scalebox{.5}{\color{blue}$1^2$}}
    (-360/9*3:-.7cm)--
    (-360/9*1:-.7cm);

    \draw[fill=lightgray] (-360/9*4:-.7cm)--
    (-360/9*5:-.7cm)--
    node[above,xshift=-.15cm,yshift=.1cm]{\scalebox{.5}{\color{blue}$1^1$}}
    (-360/9*6:-.7cm)--
    (-360/9*4:-.7cm);

    \draw[fill=lightgray] (-360/9*7:-.7cm)--
    (-360/9*8:-.7cm)--
    node[above,xshift=.25cm,yshift=-.1cm]{\scalebox{.5}{\color{blue}$1^3$}}
    (-360/9*9:-.7cm)--
    (-360/9*7:-.7cm);
    
    \draw[fill=lightgray] (-360/18*3:-2cm)--
    (-360/18*4:-2cm)--
    node[below]{\scalebox{.5}{\color{blue}$3^2$}}
    (-360/18*5:-2cm)--
    (-360/18*3:-2cm);
    
    \draw[fill=lightgray] (-360/18*9:-2cm)--
    (-360/18*10:-2cm)--
    node[above,xshift=-.1cm,yshift=.15cm]{\scalebox{.5}{\color{blue}$3^1$}}
    (-360/18*11:-2cm)--
    (-360/18*9:-2cm);

    \draw[fill=lightgray] (-360/18*15:-2cm)--
    (-360/18*16:-2cm)--
    node[above,xshift=.15cm]{\scalebox{.5}{\color{blue}$3^3$}}
    (-360/18*17:-2cm)--
    (-360/18*15:-2cm);

    \draw[fill=lightgray] (-360/18*1:-2cm)--
    (-360/18*2:-2cm)--
    node[below,xshift=-.5cm]{\scalebox{.5}{\color{blue}$2^2$}}
    (-360/18*6:-2cm)--
    (-360/18*1:-2cm);

    \draw[fill=lightgray] (-360/18*7:-2cm)--
    (-360/18*8:-2cm)--
    node[below,yshift=.5cm]{\scalebox{.5}{\color{blue}$2^1$}}
    (-360/18*12:-2cm)--
    (-360/18*7:-2cm);

    \draw[fill=lightgray] (-360/18*13:-2cm)--
    (-360/18*14:-2cm)--
    node[above,xshift=.5cm,yshift=-.5cm]{\scalebox{.5}{\color{blue}$2^3$}}
    (-360/18*18:-2cm)--
    (-360/18*13:-2cm);
    
    \foreach \a/\n in {1^1/1,2^1/2,3^1/3,4^1/4,5^1/5,6^1/6,1^3/7,2^3/8,3^3/9,4^3/10,5^3/11,6^3/12,1^2/13,2^2/14,3^2/15,4^2/16,5^2/17,6^2/18}
{
\draw[fill](-360/18*\n:-2cm)circle(1.2pt)
node[anchor=-\n*360/18] (\a) {$\a$};
}
   \foreach \a/\n in {9^1/1,8^1/2,7^1/3,9^3/4,8^3/5,7^3/6,9^2/7,8^2/8,7^2/9}
{
\draw[fill](-360/9*\n:-.7cm)circle(1.2pt)
node[anchor=-\n*360/9] (\a) {$\a$};
}
\end{tikzpicture}
            \caption{}
            \label{subfig:3mknd_b}
        \end{subfigure}\hfill
        \begin{subfigure}[h]{0.3\textwidth}
            \centering
            \begin{tikzpicture}
        \draw[fill=lightgray] (-360/18*3:-2cm)--
(-360/18*4:-2cm)--node[above,yshift=-.1cm]{\scalebox{.5}{\color{blue}$3^2$}}
(-360/18*5:-2cm)--
(-360/18*3:-2cm);

\draw[fill=lightgray] (-360/18*9:-2cm)--
(-360/18*10:-2cm)--node[below,yshift=.25cm,xshift=.15cm]{\scalebox{.5}{\color{blue}$3^1$}}
(-360/18*11:-2cm)--
(-360/18*9:-2cm);

\draw[fill=lightgray] (-360/18*15:-2cm)--
(-360/18*16:-2cm)--node[above,yshift=-.25cm,xshift=-.1cm]{\scalebox{.5}{\color{blue}$3^3$}}
(-360/18*17:-2cm)--
(-360/18*15:-2cm);

\draw[fill=lightgray] (-360/18*1:-2cm)--
(-360/18*2:-2cm)--node[below,rotate=10]{\scalebox{.5}{\color{blue}$1^1,1^2,1^3,2^1,2^2,2^3$}}
(-360/18*6:-2cm)--
(-360/18*7:-2cm)--
(-360/18*8:-2cm)--
(-360/18*12:-2cm)--
(-360/18*13:-2cm)--
(-360/18*14:-2cm)--
(-360/18*18:-2cm)--
(-360/18*1:-2cm);
   \foreach \a/\n in {1^1/1,2^1/2,3^1/3,4^1/4,5^1/5,6^1/6,1^3/7,2^3/8,3^3/9,4^3/10,5^3/11,6^3/12,1^2/13,2^2/14,3^2/15,4^2/16,5^2/17,6^2/18}
{
\draw[fill](-360/18*\n:-2cm)circle(1.2pt)
node[anchor=-\n*360/18] (\a) {$\a$};
}

 \foreach \a/\n in {9^1/1,8^1/2,7^1/3,9^3/4,8^3/5,7^3/6,9^2/7,8^2/8,7^2/9}
{
\draw[fill](-360/9*\n:-.75cm)circle(1.2pt)
node[anchor=-\n*360/9] (\a) {$\a$};
}
\end{tikzpicture}
            \caption{}
            \label{subfig:3mknd_c}
        \end{subfigure}
        \caption{Three elements of $M(k,n,d)$ with $k=3,n=4$, and $d=3$.}
        \label{fig:3_mknd_examples}
    \end{figure}

\subsection{Proof of~\Cref{thm: k-bijection} for \texorpdfstring{$G(d,d,n)$}{G(d,d,n)}}\label{sec:fuss_ddn_bij}
With $G(d,d,n)$ acting on $V = \C^n$ with the standard coordinate functionals ${\bf x}=x_1,\ldots,x_n$, choose
\[
\Theta_{kh+1} = (x_1^{kd(n-1)+1},\ldots,x_n^{kd(n-1)+1}).
\]
By~\Cref{rmk:hsop_inf_fams}, $(\Theta_{kh+1})$ is an hsop of degree $kd(n-1)+1=kh+1$ carrying $V^*(-kh-1)$ and the map $x_i\to\theta_i$ is $W$-equivariant. The subvariety cut out by $(\Theta_{kh+1} - \bf{x})$ is
\[
V^{\Theta_{kh+1}} = \{(v_1,\ldots,v_n)\in\C^n:v_i=0 \text{ or }v_i^{kd(n-1)}=1\text{ for all }i\},
\]
and we have that $|V^{\Theta_{kh+1}}|=(kh+1)^n$. We now describe the forward isomorphism in this case
\[
M(k,n,d) \xrightarrow{f_{d,k}} V^{\Theta_{kh+1}}.
\] 

To begin, we pass our partition $\pi$ restricted to $[k(n-1)]^d$ (ignoring the inner annulus) into an infinite parenthesization of the string
\[
\ldots 1^1,\ldots,(k(n-1))^1,1^d,\ldots,(k(n-1))^d,1^{d-1},\ldots,(k(n-1))^{d-1},\ldots,1^2,\ldots,(k(n-1))^d,\ldots,
\]
and extend the definition of \defn{opener} to outer and annular blocks. We define the map $f_{d,k}$ by sending $(\pi,\ell)\in M(k,n,d)$ to $(v_1,\ldots,v_n)$ where 
\begin{itemize}
    \item $v_r=0$ if $B_0\in\pi$ is the zero block and $r\in \ell(B_0)$,
    \item $v_r=0$ if $I\in\pi$ is an inner block and $r\in \ell(I)$, or 
    \item $v_r=\lambda^j\omega^i$ if $B$ is an outer or annular block opened by $i^j$ and $r^1\in \ell(B)$,
\end{itemize}
where $\omega = e^{\frac{2\pi \sqrt{-1}}{k(n-1)}}$ and $\lambda=e^{\frac{2\pi \sqrt{-1}}{d}}$.
\begin{example}
\label{ex:fussddn_forward}
The images of the three examples in~\Cref{fig:3_mknd_examples} are as follows:
\begin{itemize}
    \item[(\ref{subfig:3mknd_a})] $v=(\lambda^3\omega^5,\lambda^3\omega^1,\lambda^3\omega^1)$,
    \item[(\ref{subfig:3mknd_b})] $v=(0,\lambda^3\omega^1,\lambda^3\omega^3)$,
    \item[(\ref{subfig:3mknd_c})] $v=(0,0,\lambda^3\omega^3)$.
\end{itemize}
\end{example}

Now, to show that this map is indeed a bijection, we construct its inverse
\[
V^{\Theta_{kh+1}} \xrightarrow{f_{d,k}^{-1}} M(k,n,d).
\] Let $z$ be the number of zeroes among the components of $v$. We consider the following three cases:

\begin{itemize}
    \item[\bf{Case 1:}] $z=0$. In this case, we know $\pi$ must have only annular and outer blocks. Let $m_r$ be the number of occurrences of $\lambda^j\omega^r$ for any $j=1,\ldots, d$. Take new multiplicities $km_r-m'_r$ where $\sum_r m'_r=k$ and $km_r-m'_r\geq 1$. Any choice of $m'_r$ will recover a noncrossing partition of type $G(d,1,k(n-1))$. Similarly to the $k=1$ case, one can check that exactly one choice remains noncrossing when adding $m'_r$ many inner vertices to each block, $B_r$, corresponding to multiplicity $m_r$. To determine which of the inner vertices to add to each $B_r$, we use condition (4) from~\Cref{def: kn-noncrossing partition}. Beginning with the unique inner vertex that is congruent to the last outer vertex of $B_r$ modulo $k$, add $m'_r$ inner vertices to $B_r$. This gives $d$ potential ways to add inner vertices to each block. 

    We may obtain the labellings in the same manner as with $G(d,1,n)$. Once we have recovered the labelling in this way, we have only one choice amongst the $d$ previous choices for which condition (4) of~\Cref{lemma:mknd_conditions} holds. 

    \begin{example}
    \label{ex:Fuss_ddn_inverse_case1}
    The example in~\Cref{subfig:3mknd_a} falls into this case. Given the vector,
    \[
    v=(\lambda^3\omega^5,\lambda^3\omega^1,\lambda^3\omega^1),
    \]
    we have nonzero multiplicities $m_1=2$ and $m_5=1$. We see that $km_1-m'_1 = 3*2-2$ and $km_2-m'_2=3*1-1$ are the only choices for $m'_1$ and $m'_2$ so that the partition remains noncrossing when two central vertices are added to the block opened by $1^j$ for any color $j$ and one central vertex to the block opened by $5^j$ for any color $j$. We recover the outer blocks with the infinite parenthesization as in the $G(d,1,n)$ case. Using condition (4) of~\Cref{def: kn-noncrossing partition}, we have $d=3$ possible partitions:
    \[
    \{\{1^1,2^1,3^1,4^1,8^j,9^j\},\{5^1,6^1,7^j\}\}
    \]
    where $j=1$, $2$, or $3$. 
    \end{example}
    
    We recover the labellings in the same way as with $G(d,d,n)$, labelling the zero block with $r^j$ for all $j=1,\ldots,d$ and $r$ such that $v_r=0$. For the nonzero blocks, we do as before by labelling $B$ with $r^1$ for any $v_r = \lambda^j\omega^i$ with $i^j$ the opener of $B$. Once we have this, using condition (4) of~\Cref{lemma:mknd_conditions}, there must exist a map $\phi$ between the set $Y$ and $[n]^d$ that has a total increase of colors divisible by $d$. This condition forces one choice of our $d$ previous choices.
    \begin{example}
        Continuing with~\Cref{ex:Fuss_ddn_inverse_case1}, we have $Y=\{i^j:i=1,4,7\text{ and }j=1,2,3\}$ and the existence of such a $\phi$ restricts our choice for the block opened by $5^1$ to contain $7^1$. Thus, by the ordering of the inner boundary, we have the partition
    \[
    \{\{1^1,2^1,3^1,4^1,8^1,9^1\},\{5^1,6^1,7^1\}\}
    \]
    recovering the noncrossing parking function in~\Cref{subfig:3mknd_a}.
    \end{example}
    
    \item[\bf{Case 2:}] $z=1$. Let $m_r$ be defined as before. In this case, we must have $d$ blocks of size $k$ among the inner vertices and the remainder of the blocks are outer. By~\Cref{Thm: tau bijection}, we may use $k$ times our multiplicities $m_r$ and our known openers to reconstruct the outer blocks using the parenthesization into the infinite string up to $k(n-1)$. From here, using condition (5) in~\Cref{def: kn-noncrossing partition}, we reconstruct the partition. We recover the labellings in the same way as with $G(d,1,n)$. For the inner blocks, there are $d$ potential choices $\{\ell_i\}_{i=1}^d$ where $\ell_i(B)=\ell_j(B)$ for all $i,j$ and outer blocks $B$. Further, 
    $\ell_i(I) = \ell_j(I)^{j-i}$ for $i<j$ and any inner block $I$. All $d$ choices satisfy conditions (1)-(3) of~\Cref{lemma:mknd_conditions}, but only one satisfies condition (4). Choose this map to induce labellings and we have then constructed the inverse of $v$.  

    \begin{example}
    The example in~\Cref{subfig:3mknd_b} falls into this case. With its image
    \[
    v=(0,\lambda^3\omega^1,\lambda^3\omega^3).
    \]
    We have nonzero multiplicities $m_1=1$ and $m_3=1$. We use $km_1=3$ and $km_3=3$ to reconstruct our outer blocks using our infinite parenthesization. From here, the visible blocks are $\{1^j,2^j,6^j\}$ for $j=1,2,3$. By condition (5) of~\Cref{def: kn-noncrossing partition}, we have that the inner blocks must end with $9^j$ for any color $j$. This gives our inner blocks and we recover the labelling as before. The inner block labelling gives $d$ choices, only one of which satisfies condition (4) of~\Cref{lemma:mknd_conditions} with the set $Y=\{i^j:i=1,4,7\text{ and }j=1,2,3\}$. We thus recover our noncrossing parking function.
    \end{example}
    
    \item[\bf{Case 3:}] $z\geq2$. Let $z$ be the number of zeroes among the coordinates of $v$ and let $m_r$ be as before. In this case, we know the zero block $B_0$ of $\pi$ exists and contains all the vertices on the inner boundary of the annulus. As in case 2, we may use $k$ times our multiplicities and the known openers to reconstruct the outer boundary. This also gives what outer vertices belong to the zero block and thus our partition. We recover the labellings as in type $G(d,1,n)$, labelling the zero block by all elements not in the labellings of outer blocks.
    \begin{example}
    The example in~\Cref{subfig:3mknd_c} falls into this case. With the vector,
    \[
    v=(0,0,\lambda^3\omega^3),
    \]
    we must have a zero block and the only nonzero multiplicity is $m_3=1$. Using our infinite parenthesization, we recover the outer block and the vertices on the outer boundary that belong to the zero block. We recover the labellings in the same way as with $G(d,1,n)$, labelling the zero block with any unused letters.
    \end{example}
\end{itemize}

We have thus shown this map is a bijection. We leave it to the reader to verify this map is $(W\times \Z_{kh})$-equivariant. Thus, we have proven~\Cref{thm: k-bijection} in this case and~\Cref{thm: k-Park is (h+1)^n} by our choice of hsop. 

\section{Future Work}\label{sec:future_work}
We plan to continue our study relating to the noncrossing parking functions in complex reflection groups. Recent work has been done in regard to the poset topology of the noncrossing parking functions for finite real reflection groups~\cite{Douv2024clusterparking}. We plan to extend this work to well--generated complex reflection groups.

First, let us briefly introduce some notions of poset topology. Let $\Gamma$ be a finite simplicial complex, identified with its face poset with inclusion its order relation. We say $\Gamma$ is \defn{pure} if its facets (the maximal faces) all have the same dimension. Let $P$ be a finite poset. If $P$ has a unique minimal element and a unique maximal element, we will write $\overline{P}$ to denote the \defn{proper part} obtained by removing said unique minimal and maximal elements. For $P$ a finite poset, the \defn{order complex}, $\Omega(P)$, of $P$ is the simplicial complex having $P$ as vertices and pairwise comparable elements as its faces. 
 
Let $W$ be a finite real reflection group (not necessarily irreducible) of rank $n$ and having Coxeter number $h$. The noncrossing parking function poset was initially introduced by P.~Edelman~\cite{edelman1980chain} as $2$-partitions and further studied in~\cite{Delcroix2022pfposet}. In regards to the topology of the noncrossing parking functions, B.~Delcroix-Oger, M.~Josuat-Verg\"es, and L.~Randazzo handle the type $A$ case using the shellability of the poset~\cite{Delcroix2022pfposet}. For general finite real reflection groups, T.~Douvropoulos and M.~Josuat-Verg\"es describe the homotopy type of the noncrossing parking function poset as a pure wedge of spheres dependent on $n$ and $h$~\cite{Douv2024clusterparking}, obtaining the following result:
\begin{theorem}[\cite{Douv2024clusterparking}]
    For $W$ a finite real reflection group, $\Omega(\overline{\Park^\NC (W)})\simeq (\mathbb{S}^{n-1})^{\vee (h-1)^n}$,
    where $\mathbb{S}^m$ is the $m$-dimensional sphere.
\end{theorem}

Similar to the work in this paper, we wish to extend the work of T.~Douvropoulos and M.~Josuat-Verg\"es in~\cite{Douv2024clusterparking} to well--generated complex reflection groups. 

\begin{conjecture}
    For $W$ a well--generated complex reflection group of rank $n$ and Coxeter number $h$,
    \[
    \Omega(\overline{\Park^\NC (W)})\simeq (\mathbb{S}^{n-1})^{\vee(h-1)^n}.
    \] 
\end{conjecture}

\section{Acknowledgments}
The author would like to thank Nathan Williams for his guidance and many helpful discussions regarding this work.  The author was partially supported by the National
Science Foundation under Grant No.~DMS-2246877.

\bibliographystyle{amsalpha}
\bibliography{lib}

\end{document}